# THE LIL FOR CANONICAL $U$-STATISTICS

By Radosław Adamczak[1] and Rafał Latała[2]

*Polish Academy of Sciences and Warsaw University*

We give necessary and sufficient conditions for the (bounded) law of the iterated logarithm for canonical $U$-statistics of arbitrary order $d$, extending the previously known results for $d = 2$. The nasc's are expressed as growth conditions on a parameterized family of norms associated with the $U$-statistics kernel.

**1. Introduction.** $U$-statistics [i.e. statistics being averages of a measurable kernel $h(x_1, \ldots, x_d)$ over an i.i.d. sample $X_1, X_2, \ldots, X_n$] were introduced by Hoeffding [11] and Halmos [9] in the 1940s and since then have become an important tool in asymptotic statistics, appearing for instance as unbiased estimators or higher-order terms in expansions of smooth statistics. Their relevance stems mainly from the fact that they share many basic properties with sums of i.i.d. random variables. Already in the 1960s Hoeffding proved that $\mathbb{E}|h| < \infty$ is a sufficient condition for a $U$-statistic to satisfy the SLLN [12], the CLT under the finiteness of the second moment of the kernel (and complete degeneracy—a technical assumption which will be explained in the sequel) was obtained by Rubin and Vitale in 1980 [18], finally the LIL (under the same hypothesis) was proved by Arcones and Giné in 1995 [2]. All the abovementioned results are occurrences of a general phenomenon, manifesting itself in the fact that the necessary and sufficient conditions for the classical triple of limit theorems for sums of i.i.d. random variables (SLLN, CLT or LIL) are sufficient for analogous limit theorems for $U$-statistics. It may be, therefore, somewhat surprising (and as a matter of fact remained for some time unnoticed) that with the exception of the CLT, these conditions fail to be necessary.

Recently we have witnessed a rapid development in the asymptotic theory of $U$-statistics, following the discovery of the so-called decoupling technique

Received April 2006; revised May 2007.
[1]Supported in part by MEiN Grant 2 PO3A 019 30.
[2]Supported in part by MEiN Grant 1 PO3A 012 29.

(see [3] and the references therein), which allows one to treat $U$-statistics as sums of conditionally independent random variables. In particular, the sufficient conditions for the CLT given by Rubin and Vitale were proven to be also necessary (Giné and Zinn [7]). Also the necessary and sufficient conditions for the SLLN were found ([19] for $d=2$, [15] for general $d$). In 1999 Giné et al. [8] obtained necessary and sufficient conditions for the law of the iterated logarithm for $U$-statistics of order 2. The conditions they gave turned out to be less restrictive and more subtle than just the square integrability of the kernel (as indicated already by Giné and Zhang [5]). Completing the picture requires finding the nasc's for the LIL in the general case and identifying the limit set in the LIL (which in general is unknown even for $d=2$).

In this paper, we address the first of these questions, namely we give the nasc's on a kernel $h(x_1,\dots,x_d)$ to satisfy the (bounded) law of the iterated logarithm. In particular we prove that a conjecture stated in [8] is false.

**2. Notation.** For an integer $d$, let $(X_i)_{i\in\mathbb{N}}$, $(X_i^{(k)})_{i\in\mathbb{N},1\le k\le d}$ be i.i.d. random variables with values in a Polish space $\Sigma$, equipped with the Borel $\sigma$-field $\mathcal{F}$. Consider moreover a measurable function $h:\Sigma^d\to\mathbb{R}$.

To shorten the notation, we will use the following convention. For $\mathrm{i}=(i_1,\dots,i_d)\in\{1,\dots,n\}^d$ we will write $\mathbf{X}_\mathrm{i}$ (resp. $\mathbf{X}_\mathrm{i}^{\mathrm{dec}}$) for $(X_{i_1},\dots,X_{i_d})$ (resp. $(X_{i_1}^{(1)},\dots,X_{i_d}^{(d)})$) and $\epsilon_\mathrm{i}$ (resp. $\epsilon_\mathrm{i}^{\mathrm{dec}}$) for the product $\varepsilon_{i_1}\cdot\dots\cdot\varepsilon_{i_d}$ (resp. $\varepsilon_{i_1}^{(1)}\cdot\dots\cdot\varepsilon_{i_d}^{(d)}$), the notation being thus slightly inconsistent, which however should not lead to a misunderstanding. The $U$-statistics will, therefore, be denoted

$$\sum_{\mathrm{i}\in I_n^d} h(\mathbf{X}_\mathrm{i}) \qquad \text{(an undecoupled $U$-statistic)}$$

$$\sum_{|\mathrm{i}|\le n} h(\mathbf{X}_\mathrm{i}^{\mathrm{dec}}) \qquad \text{(a decoupled $U$-statistic)}$$

$$\sum_{\mathrm{i}\in I_n^d} \epsilon_\mathrm{i} h(\mathbf{X}_\mathrm{i}) \qquad \text{(an undecoupled randomized $U$-statistic)}$$

$$\sum_{|\mathrm{i}|\le n} \epsilon_\mathrm{i}^{\mathrm{dec}} h(\mathbf{X}_\mathrm{i}^{\mathrm{dec}}) \qquad \text{(a decoupled randomized $U$-statistic)},$$

where

$$|\mathrm{i}| = \max_{k=1,\dots,d} i_k,$$

$$I_n^d = \{\mathrm{i}: |\mathrm{i}|\le n, i_j\ne i_k \text{ for } j\ne k\}.$$

Since in this notation $\{1,\dots,d\}=I_d^1$ we will write

$$I_d=\{1,2,\dots,d\}.$$

We will also occasionally write $X$ for $(X_1,\dots,X_d)$ and for $I\subseteq I_d$, $X_I=(X_i)_{i\in I}$. Sometimes we will write simply $h$ instead of $h(X)$.

Throughout the article we will write $K, L_d, L$ to denote constants depending only on the function $h$, only on $d$ and universal constants, respectively. In all those cases the values of a constant may differ at each occurrence.

To avoid technical problems with small values of $h$ let us also define $\mathrm{LL}\,x=\log\log(x\vee e^e)$.

Let us also introduce some notation for conditional expectation. For $j\in I_d$, by $\mathbb{E}_j$ we will denote expectation with respect to $(X_i^{(j)})_i$, $((\varepsilon_i^{(j)},X_i^{(j)}))_i$ or $X_j$ (depending on the context). Similarly, for $I\subseteq I_d$, we will denote by $\mathbb{E}_I$, integration with respect to $(X_i^{(j)})_{j\in I,i}$, $((\varepsilon_i^{(j)},X_i^{(j)}))_{j\in I,i}$ or $(X_i)_{i\in I}$. Although at first this notation may seem slightly ambiguous, it turns out to be quite natural at specific instances and should not lead to misunderstanding.

In the article we will consider mainly *canonical* (or *completely degenerate*) kernels, that is kernels $h$, such that for all $j\in I_d$,

$$\mathbb{E}_j h(X_1,\dots,X_d)=0 \qquad \text{a.s.}$$

**3. The main result.** Let us now introduce the quantities, that the necessary and sufficient conditions for the LIL will be expressed in.

DEFINITION 1. For a finite set $I$, let $\mathcal{P}_I$ denote the family of all partitions of $I$ into disjoint, nonempty sets and for a partition $\mathcal{J}\in\mathcal{P}_I$ let $\deg\mathcal{J}$ be the number of elements of $\mathcal{J}$. For a kernel $h:\Sigma^d\to\mathbb{R}$, a partition $\mathcal{J}=\{J_1,\dots,J_k\}\in\mathcal{P}_{I_d}$ and a nonnegative number $u$, define

$$\begin{aligned}\|h\|_{\mathcal{J},u}&=\|h(X)\|_{\mathcal{J},u}\\&=\sup\biggl\{\mathbb{E}\biggl[h(X)\prod_{i=1}^k f_i(X_{J_i})\biggr]:\|f_i(X_{J_i})\|_2\le 1,\\&\qquad\qquad \|f_i(X_{J_i})\|_\infty\le u, i=1,\dots,k\biggr\}.\end{aligned}$$

EXAMPLE. For $d=3$, the above definition gives

$$\begin{aligned}\|h(X_1,X_2,X_3)\|_{\{1,2,3\},u}&=\sup\{\mathbb{E}h(X_1,X_2,X_3)f(X_1,X_2,X_3):\\&\qquad \mathbb{E}f(X_1,X_2,X_3)^2\le 1,\|f\|_\infty\le u\},\\ \|h(X_1,X_2,X_3)\|_{\{1,2\}\{3\},u}&=\sup\{\mathbb{E}h(X_1,X_2,X_3)f(X_1,X_2)g(X_3):\\&\qquad \mathbb{E}f(X_1,X_2)^2,\mathbb{E}g(X_3)^2\le 1,\\&\qquad \|f\|_\infty,\|g\|_\infty\le u\},\end{aligned}$$

$$\begin{aligned}\|h(X_1,X_2,X_3)\|_{\{1\}\{2\}\{3\},u} = \sup\{&\mathbb{E}h(X_1,X_2,X_3)f(X_1)g(X_2)k(X_3)\colon\\ &\mathbb{E}f(X_1)^2,\mathbb{E}g(X_2)^2,\mathbb{E}k(X_3)^2\le 1,\\ &\|f\|_\infty,\|g\|_\infty,\|k\|_\infty\le u\}.\end{aligned}$$

Although at first approach the $\|\cdot\|_{\mathcal{J},u}$ norms may seem quite unusual, they resemble both the quantities appearing in tail estimates for canonical $U$-statistics and in tail estimates for Rademacher chaoses (see Sections 4.2 and 4.3 below) and they indeed play an important role in necessary and sufficient conditions for the LIL, as can be seen in our main result, which is

THEOREM 1. *For any symmetric $h\colon\Sigma^d\to\mathbb{R}$, the law of the iterated logarithm*

$$\limsup_{n\to\infty}\frac{1}{(n\log\log n)^{d/2}}\left|\sum_{\mathbf{i}\in I_n^d}h(\mathbf{X_i})\right|<\infty \qquad a.s.$$

*holds if and only if $h$ is completely degenerate for the law of $X_1$ and for all $\mathcal{J}\in\mathcal{P}_{I_d}$,*

$$\limsup_{u\to\infty}\frac{1}{(\log\log u)^{(d-\deg\mathcal{J})/2}}\|h\|_{\mathcal{J},u}<\infty.$$

*(Recall that according to Definition 1, $\deg\mathcal{J}$ denotes the number of elements of $\mathcal{J}$.)*

REMARK. Obviously, although formally in the above theorem one considers all the partitions $\mathcal{J}$, due to symmetry of the kernel and equidistribution of the variables $X_1,\ldots,X_d$, many of them give the same value of $\|h\|_{\mathcal{J},u}$. For instance for $d=3$ we have $\|h\|_{\{1\}\{2,3\},u}=\|h\|_{\{2\}\{1,3\},u}=\|h\|_{\{3\}\{1,2\},u}$ (note that we suppressed the outer brackets in the lower index and wrote e.g. $\|h\|_{\{2\}\{1,3\},u}$ instead of $\|h\|_{\{\{2\}\{1,3\}\},u}$. We will do so whenever there is no risk of confusion also with similar norms, which will be introduced in Sections 4.2 and 4.3).

## 4. Preliminaries. Basic definitions and tools.

4.1. *Hoeffding's decomposition.* We will now describe a decomposition of a $U$-statistic with mean zero kernel into a sum of completely degenerate $U$-statistics, introduced in [11], which is one of the basic tools in the analysis of $U$-statistics. Recall that we are working with a fixed sequence $(X_i)_{i\in\mathbb{N}}$ of i.i.d. $\Sigma$-valued random variables. Then the classical definition of Hoeffding's projections is as follows.

DEFINITION 2. For an integrable kernel $h:\Sigma^d\to\mathbb{R}$ and $k=0,1,\dots,d$, define $\pi_k h:\Sigma^k\to\mathbb{R}$ with the formula

$$\pi_k h(x_1,\dots,x_k)=(\delta_{x_1}-\mathbf{P})\times(\delta_{x_2}-\mathbf{P})\times\cdots\times(\delta_{x_k}-\mathbf{P})\times\mathbf{P}^{d-k}h,$$

where $\mathbf{P}$ is the law of $X_1$.

In particular $\pi_0 h=\mathbb{E}h, \pi_1 h(x_1)=\mathbb{E}_{\{2,\dots,d\}}h(x_1,X_2,\dots,X_d)-\mathbb{E}h$.

We will however need to extend this definition (for $k=d$) to $U$-statistics based not necessarily on an i.i.d. sequence. Let us thus introduce the following definition

DEFINITION 3. Let $h:\Sigma_1\times\cdots\times\Sigma_d\to\mathbb{R}$ be a measurable function. Consider independent sequences $(X_j^{(1)})_j,\dots,(X_j^{(d)})_j$ of i.i.d. random variables with values in $\Sigma_1,\dots,\Sigma_d$ respectively, such that $\mathbb{E}|h(X_1^{(1)},\dots,X_1^{(d)})|<\infty$. Define $\pi_d h:\Sigma_1\times\cdots\times\Sigma_d\to\mathbb{R}$ with the formula

$$\pi_d h(x_1,\dots,x_d)=(\delta_{x_1}-\mathbf{P}_{X_1^{(1)}})\times\cdots\times(\delta_{x_d}-\mathbf{P}_{X_1^{(d)}})h,$$

where $\mathbf{P}_{X_1^{(i)}}$ is the law of $X_1^{(i)}$.

Obviously for $\Sigma_1=\cdots=\Sigma_d$ and $(X_i^{(j)})_{i\in\mathbb{N}}$—independent copies of $(X_i)_{i\in\mathbb{N}}$, the above definitions of $\pi_d h$ are equivalent.

It is easy to check that for $k\geq 1$, $\pi_k h$ is canonical for the law of $X_1$ (note also that $\pi_0 h=\mathbb{E}h$).

In the sequel we will need the following comparison of moments for $U$-statistics:

LEMMA 1. *Consider an arbitrary family of integrable kernels $h_{\mathrm{i}}:\Sigma_1\times\cdots\times\Sigma_d\to\mathbb{R}$, $|\mathrm{i}|\leq n$. For any $p\geq 1$ we have*

$$\left\|\sum_{|\mathrm{i}|\leq n}\pi_d h_{\mathrm{i}}(\mathbf{X}_{\mathrm{i}}^{\mathrm{dec}})\right\|_p\leq 2^d\left\|\sum_{|\mathrm{i}|\leq n}\epsilon_{\mathrm{i}}^{\mathrm{dec}}h_{\mathrm{i}}(\mathbf{X}_{\mathrm{i}}^{\mathrm{dec}})\right\|_p.$$

PROOF. For $d=1$, the statement of the lemma is the classical symmetrization inequality for sums of independent random variables. Now, we use induction with respect to $d$. To simplify the notation let $\bar{\pi}_{d-1}h_{\mathrm{i}}$ denote the proper Hoeffding's projection of $h_{\mathrm{i}}$ treated as a function of $x_2,\dots,x_d$, with the first coordinate fixed, that is

$$\bar{\pi}_{d-1}h_{\mathrm{i}}(x)=\delta_{x_1}\times(\delta_{x_2}-\mathbf{P}_{X_1^{(2)}})\times\cdots\times(\delta_{x_d}-\mathbf{P}_{X_1^{(d)}})h_{\mathrm{i}}.$$

Assume now that the lemma is true for all kernels of degree smaller than $d$. Consider $(\tilde{X}_i^{(k)})_{i\in\mathbb{N},k\leq d}$, an independent copy of $(X_i^{(k)})_{i\in\mathbb{N},k\leq d}$ and denote by $\tilde{\mathbb{E}}_1$ integration with respect to $\tilde{X}^{(1)}$. Then, the complete degeneracy of $\pi_d h_{\mathrm{i}}$ and Jensen's inequality yield

$$\begin{aligned}
\mathbb{E}_1&\Big|\sum_{|\mathrm{i}|\le n}\pi_d h_{\mathrm{i}}(\mathbf{X}_{\mathrm{i}}^{\mathrm{dec}})\Big|^p\\
&=\mathbb{E}_1\Big|\sum_{|\mathrm{i}|\le n}(\pi_d h_{\mathrm{i}}(X_{i_1}^{(1)},\dots,X_{i_d}^{(d)})-\tilde{\mathbb{E}}_1\pi_d h_{\mathrm{i}}(\tilde{X}_{i_1}^{(1)},X_{i_2}^{(2)},\dots,X_{i_d}^{(d)}))\Big|^p\\
&\le\mathbb{E}_1\tilde{\mathbb{E}}_1\Big|\sum_{|\mathrm{i}|\le n}(\pi_d h_{\mathrm{i}}(X_{i_1}^{(1)},\dots,X_{i_d}^{(d)})-\pi_d h_{\mathrm{i}}(\tilde{X}_{i_1}^{(1)},X_{i_2}^{(2)},\dots,X_{i_d}^{(d)}))\Big|^p\\
&=\mathbb{E}_1\tilde{\mathbb{E}}_1\Big|\sum_{|\mathrm{i}|\le n}\varepsilon_{i_1}^{(1)}(\pi_d h_{\mathrm{i}}(X_{i_1}^{(1)},\dots,X_{i_d}^{(d)})-\pi_d h_{\mathrm{i}}(\tilde{X}_{i_1}^{(1)},X_{i_2}^{(2)},\dots,X_{i_d}^{(d)}))\Big|^p\\
&=\mathbb{E}_1\tilde{\mathbb{E}}_1\Big|\sum_{|\mathrm{i}|\le n}\varepsilon_{i_1}^{(1)}(\bar{\pi}_{d-1} h_{\mathrm{i}}(X_{i_1}^{(1)},\dots,X_{i_d}^{(d)})\\
&\qquad\qquad\qquad-\bar{\pi}_{d-1} h_{\mathrm{i}}(\tilde{X}_{i_1}^{(1)},X_{i_2}^{(2)},\dots,X_{i_d}^{(d)}))\Big|^p,
\end{aligned}$$

so, using the triangle inequality, we obtain

$$\Big\|\sum_{|\mathrm{i}|\le n}\pi_d h_{\mathrm{i}}(\mathbf{X}_{\mathrm{i}}^{\mathrm{dec}})\Big\|_p\le 2\Big\|\sum_{|\mathrm{i}|\le n}\varepsilon_{i_1}^{(1)}\bar{\pi}_{d-1} h_{\mathrm{i}}(\mathbf{X}_{\mathrm{i}}^{\mathrm{dec}})\Big\|_p.$$

Now, the Fubini theorem, together with the induction assumption applied to the family of kernels $\tilde{h}_{(i_2,\dots,i_d)}(x_2,\dots,x_d)=\sum_{i_1\le n}\varepsilon_{i_1}^{(1)}h_{\mathrm{i}}(X_{i_1}^{(1)},x_2,\dots,x_d)$ for fixed values of $X^{(1)},\varepsilon^{(1)}$, proves the lemma. □

We will also use the classical theorem due to Hoeffding, giving a decomposition of a $U$-statistic into sum of uncorrelated, canonical $U$-statistics of different orders, mentioned at the beginning of this paragraph.

LEMMA 2 (Hoeffding's decomposition; see, e.g. [3], page 137). *For $h:\Sigma^d\to\mathbb{R}$, symmetric in its entries denote*

$$U_n(h)=\frac{(n-d)!}{n!}\sum_{\mathrm{i}\in I_n^d}h(\mathbf{X}_{\mathrm{i}}).$$

*Then*

$$U_n(h)=\sum_{k=0}^d\binom{d}{k}U_n(\pi_k h).$$

4.2. *Moment and tail estimates for canonical $U$-statistics.* We will now present a version of sharp moment estimates for canonical $U$-statistics, proved in [1] (actually as we will not need these results in the whole generality, we will state only a simplified corollary, adapted to our purposes, which follows immediately from Theorem 6 there).

First let us introduce some quantities, which will appear in the moment estimates.

DEFINITION 4. For any canonical kernel $h:\Sigma^d \to \mathbb{R}$ and each $\mathcal{J} = \{J_1, \ldots, J_k\} \in \mathcal{P}_{I_d}$ define the norm

$$\|h\|_{\mathcal{J}} := \|h\|_{\mathcal{J},\infty} = \sup\left\{\mathbb{E}\left[h(X)\prod_{i=1}^{k} f_i(X_{J_i})\right] : \mathbb{E}f_i(X_{J_i})^2 \le 1, i = 1, \ldots, k\right\}.$$

Thus $\|h\|_{\mathcal{J}}$ is the norm of $h$ viewed as a $k$-linear functional acting on the space $L^2(X_{J_1}) \times \cdots \times L^2(X_{J_k})$, where $L^2(X_{J_i})$ is the space of all square integrable random variables, measurable with respect to $\sigma(X_{J_i})$, the $\sigma$-field generated by $X_{J_i}$. In particular $\|h\|_{I_d} = (\mathbb{E}h^2)^{1/2}$ and $\|h\|_{\{1\}\ldots\{d\}}$ is the norm of $h$ seen as a kernel of a $d$-linear functional.

We have the following (cf. [1], Theorem 6)

THEOREM 2. *There exist constants $L_d$, such that for all canonical kernels $h:\Sigma^d \to \mathbb{R}$ and $p \ge 2$,*

$$\mathbb{E}\left|\sum_{|\mathrm{i}|\le n} h(\mathbf{X}_{\mathrm{i}}^{\mathrm{dec}})\right|^p \le L_d^p\Bigg[n^{dp/2}\sum_{\mathcal{J}\in\mathcal{P}_{I_d}} p^{p\deg(J)/2}\|h\|_{\mathcal{J}}^p + \sum_{I\subsetneq I_d} n^{p\#I/2}p^{p(d+\#I^c)/2}\mathbb{E}_{I^c}\max_{\mathrm{i}_{I^c}}(\mathbb{E}_I h(\mathbf{X}_{\mathrm{i}}^{\mathrm{dec}})^2)^{p/2}\Bigg].$$

REMARK. Note that $(\mathbb{E}_I h(\mathbf{X}_{\mathrm{i}}^{\mathrm{dec}})^2)^{p/2}$ depends only on $X_{\mathrm{i}_{I^c}}$, so the expression $\max_{\mathrm{i}_{I^c}}(\mathbb{E}_I h(\mathbf{X}_{\mathrm{i}}^{\mathrm{dec}})^2)^{p/2}$ in the above inequality is well defined.

Theorem 2 implies the following theorem.

THEOREM 3. *There exist constants $L_d$, such that for all bounded, canonical kernels $h:\Sigma^d \to \mathbb{R}$ and $t \ge 0$,*

$$\mathbb{P}\left(\left|\sum_{|\mathrm{i}|\le n} h(\mathbf{X}_{\mathrm{i}}^{\mathrm{dec}})\right| \ge t\right) \le L_d\exp\Bigg[-\frac{1}{L_d}\Bigg(\min_{\mathcal{J}\in\mathcal{P}_{I_d}}\left(\frac{t}{n^{d/2}\|h\|_{\mathcal{J}}}\right)^{2/\deg(\mathcal{J})}$$

$$\wedge \min_{I \subsetneq I_d} \left( \frac{t}{n^{\# I/2} \| (\mathbb{E}_I h^2)^{1/2} \|_\infty} \right)^{2/(d+\# I^c)} \Bigg) \Bigg].$$

REMARK. We would like to stress that Theorem 3 has been obtained from Theorem 2 by means of the Chebyshev inequality only. Therefore, the same tail estimates hold for random variables whose moments are dominated by moments of corresponding $U$-statistics, which together with Lemma 1 yields the following.

THEOREM 4. *There exist constants $L_d$, such that for all bounded kernels $h : \Sigma^d \to \mathbb{R}$ and all $t \geq 0$,*

$$\mathbb{P}\left( \left| \sum_{|\mathrm{i}| \leq n} \pi_d h(\mathbf{X}_{\mathrm{i}}^{\mathrm{dec}}) \right| \geq t \right)$$
$$\leq L_d \exp\Bigg[ -\frac{1}{L_d} \Bigg( \min_{\mathcal{J} \in \mathcal{P}_{I_d}} \left( \frac{t}{n^{d/2} \|h\|_{\mathcal{J}}} \right)^{2/\deg(\mathcal{J})}$$
$$\wedge \min_{I \subsetneq I_d} \left( \frac{t}{n^{\# I/2} \| (\mathbb{E}_I h^2)^{1/2} \|_\infty} \right)^{2/(d+\# I^c)} \Bigg) \Bigg].$$

4.3. *Moment and tail estimates for Rademacher chaoses.*

LEMMA 3. *Let $(a_{\mathrm{i}})_{\mathrm{i} \in I_n^d}$ be a $d$-indexed array of real numbers. Let us consider a random variable*

$$S := \left| \sum_{|\mathrm{i}| \leq n} a_{\mathrm{i}} \prod_{k=1}^{d} \varepsilon_{i_k}^{(k)} \right| = \left| \sum_{|\mathrm{i}| \leq n} a_{\mathrm{i}} \epsilon_{\mathrm{i}}^{\mathrm{dec}} \right|.$$

*Moreover, for any partition $\mathcal{J} = \{J_1, \ldots, J_m\} \in \mathcal{P}_{I_d}$ let us define*

$$\|(a_{\mathrm{i}})\|_{\mathcal{J},p}^* := \sup\Bigg\{ \left| \sum_{|\mathrm{i}| \leq n} a_{\mathrm{i}} \prod_{k=1}^{m} \alpha_{\mathrm{i}_{J_k}}^{(k)} \right| : \sum_{\mathrm{i}_{J_k}} (\alpha_{\mathrm{i}_{J_k}}^{(k)})^2 \leq p,$$
$$\forall_{i_{\max J_k} \in I_n} \sum_{\mathrm{i}_{\diamond J_k}} (\alpha_{\mathrm{i}_{J_k}}^{(k)})^2 \leq 1, k = 1, \ldots, m \Bigg\},$$

*where $\diamond J_k = J_k \backslash \{\max J_k\}$ (here $\sum_{\mathrm{i}_\varnothing} a_{\mathrm{i}} = a_{\mathrm{i}}$). Then, for all $p \geq 1$,*

$$\|S\|_p \geq \frac{1}{L_d} \sum_{\mathcal{J} \in \mathcal{P}_{I_d}} \|(a_{\mathrm{i}})\|_{\mathcal{J},p}^*.$$

*In particular for some constant $c_d$,*

$$\mathbb{P}\left( S \geq c_d \sum_{\mathcal{J} \in \mathcal{P}_{I_d}} \|(a_{\mathrm{i}})\|_{\mathcal{J},p}^* \right) \geq c_d \wedge e^{-p}.$$

PROOF. We will use induction with respect to $d$. For $d=1$ the inequalities of the lemma have been proved in [10], for $d=2$ in [14] (as a part of much sharper two-sided inequalities). Let us thus assume that the moment estimate holds for chaoses of order smaller than $d\geq 3$.

First, consider the partition $\mathcal{J}=\{I_d\}$. We have

$$\begin{aligned}\mathbb{E}S^p &= \mathbb{E}_d\mathbb{E}_{I_{d-1}}\left|\sum_{i_d}\varepsilon_{i_d}^{(d)}\left|\sum_{\mathrm{i}_{I_{d-1}}} a_{\mathrm{i}}\prod_{k=1}^{d-1}\varepsilon_{i_k}^{(k)}\right|\right|^p \\ &\geq \mathbb{E}_d\left|\sum_{i_d}\varepsilon_{i_d}^{(d)}\mathbb{E}_{I_{d-1}}\left|\sum_{\mathrm{i}_{I_{d-1}}} a_{\mathrm{i}}\prod_{k=1}^{d-1}\varepsilon_{i_k}^{(k)}\right|\right|^p \\ &\geq \frac{1}{\tilde{L}_d^p}\mathbb{E}_d\left|\sum_{i_d}\varepsilon_{i_d}^{(d)}\left(\sum_{\mathrm{i}_{I_{d-1}}} a_{\mathrm{i}}^2\right)^{1/2}\right|^p \\ &\geq \frac{1}{\tilde{L}_d^p L_1^p}\sup\left\{\sum_{i_d}\alpha_{i_d}\left(\sum_{\mathrm{i}_{I_{d-1}}} a_{\mathrm{i}}^2\right)^{1/2} : \sum_{i_d}\alpha_{i_d}^2\leq p, |\alpha_i|\leq 1\right\}^p,\end{aligned}$$

where the first inequality follows from Jensen's inequality, the second one from hypercontractivity of Rademacher chaos (see [3], Theorem 3.2.5, page 115) and the contraction principle for Rademacher averages (see, for instance, [16], Theorem 4.4, page 95), whereas the third one follows from the induction assumption.

It remains to show that

$$\sup\left\{\sum_{i_d}\alpha_{i_d}\left(\sum_{\mathrm{i}_{I_{d-1}}} a_{\mathrm{i}}^2\right)^{1/2} : \sum_{i_d}\alpha_{i_d}^2\leq p, |\alpha_i|\leq 1\right\}\geq \|(a_{\mathrm{i}})\|^*_{\{I_d\},p}.$$

Let thus $(\gamma_{\mathrm{i}})$ be a $d$-indexed matrix, such that $\sum_{\mathrm{i}}\gamma_{\mathrm{i}}^2\leq p$, $\sum_{\mathrm{i}_{I_{d-1}}}\gamma_{\mathrm{i}}^2\leq 1$ for all $i_d$. Then

$$\begin{aligned}\left|\sum_{\mathrm{i}}\gamma_{\mathrm{i}}a_{\mathrm{i}}\right| &\leq \sum_{\mathrm{i}}|\gamma_{\mathrm{i}}||a_{\mathrm{i}}| \leq \sum_{i_d}\left(\sum_{\mathrm{i}_{I_{d-1}}}\gamma_{\mathrm{i}}^2\right)^{1/2}\left(\sum_{\mathrm{i}_{I_{d-1}}} a_{\mathrm{i}}^2\right)^{1/2} \\ &\leq \sup\left\{\sum_{i_d}\alpha_{i_d}\left(\sum_{\mathrm{i}_{I_{d-1}}} a_{\mathrm{i}}^2\right)^{1/2} : \sum_{i_d}\alpha_{i_d}^2\leq p, |\alpha_i|\leq 1\right\}.\end{aligned}$$

Let now $\mathcal{J}=\{J_1,\ldots,J_m\}$, $m\geq 2$. We have

$$\|S\|_p \geq \frac{1}{L_{d-\#J_1}}\left(\mathbb{E}_{J_1}\left(\left\|\left(\sum_{\mathrm{i}_{J_1}} a_{\mathrm{i}}\prod_{k\in J_1}\varepsilon_{i_k}^{(k)}\right)_{\mathrm{i}_{I_d\setminus J_1}}\right\|^*_{\mathcal{J}\setminus\{J_1\},p}\right)^p\right)^{1/p}$$

$$\geq \frac{1}{L_{d-\#J_1}L_{\#J_1}} \|(a_{\mathbf{i}})\|^*_{\mathcal{J},p},$$

by the induction assumption and Jensen's inequality. □

4.4. *Basic consequences of the integrability condition.* Now we would like to present some basic facts, following from the integrability condition $\mathbb{E}(h^2 \wedge u) = \mathcal{O}((\log\log u)^{d-1})$, which is necessary for the LIL for $U$-statistics of order $d$, as proved by Giné and Zhang [5]; cf. Lemma 7 below.

LEMMA 4. *If $\mathbb{E}(h^2 \wedge u) = \mathcal{O}((\log\log u)^{d-1})$ then for $I \subseteq I_d$, $I \neq \varnothing, I_d$ and $a > 0$,*

$$\sum_{l=0}^{\infty}\sum_{n=3}^{\infty} 2^{l+\#I^c n}\mathbb{P}_{I^c}(\mathbb{E}_I(h^2 \wedge 2^{an}) \geq 2^{2l+\#I^c n}\log^d n) < \infty.$$

*As a consequence, for $k \geq 0$,*

$$\sum_n 2^{\#I^c n}(\log n)^{-k}\mathbb{P}_{I^c}(\mathbb{E}_I(h^2 \wedge 2^{2nd}) \geq 2^{\#I^c n}(\log n)^{d-k}) < \infty.$$

PROOF. For fixed $l$ and $k$ we have

$$\begin{aligned}
&\sum_{2^k < \log n \leq 2^{k+1}} 2^{l+\#I^c n}\mathbb{P}_{I^c}(\mathbb{E}_I(h^2 \wedge 2^{an}) \geq 2^{2l+\#I^c n}\log^d n)\\
&\quad\leq \sum_{2^k < \log n \leq 2^{k+1}} 2^{l+\#I^c n}\mathbb{P}_{I^c}(\mathbb{E}_I(h^2 \wedge 2^{ae^{2^{k+1}}}) \geq 2^{2l+\#I^c n+dk})\\
&\quad\leq 2^l\mathbb{E}_{I^c}\sum_n 2^{\#I^c n}\mathbf{1}_{\{\mathbb{E}_I(h^2 \wedge 2^{ae^{2^{k+1}}}) \geq 2^{2l+\#I^c n+dk}\}}\\
&\quad\leq 2^l\mathbb{E}_{I^c}\left[2\frac{\mathbb{E}_I(h^2 \wedge 2^{ae^{2^{k+1}}})}{2^{2l+dk}}\right] \leq 2^{-l}K\frac{(\log ae^{2^{k+1}})^{d-1}}{2^{dk}}\\
&\quad\leq 2^{-l}\tilde{K}\left(\frac{\log^{d-1} a}{2^{dk}} + 2^{-k}\right),
\end{aligned}$$

with $\tilde{K}$ depending only on $h$ (recall the convention explained in Section 2), which proves the first part of the lemma. To obtain the other inequality, it is enough to make an approximate change of variable $2^{\#I^c m}(\log m)^{-k} \simeq 2^{\#I^c n}$ and use the convergence of the inner sum for $l = 0$ in the first inequality, for $a > 2d$. □

LEMMA 5. *If $\mathbb{E}(h^2 \wedge u) = \mathcal{O}((\log\log u)^{\beta})$ then*

$$\mathbb{E}|h|\mathbf{1}_{\{|h| \geq s\}} = \mathcal{O}\left(\frac{(\log\log s)^{\beta}}{s}\right).$$

PROOF. Indeed, since $\mathbb{P}(|h| \geq u) \leq K \frac{(\log\log u)^\beta}{u^2}$, we have for sufficiently large $s$,

$$\mathbb{E}|h|\mathbf{1}_{\{|h|\geq s\}} = \sum_{k=0}^{\infty} \mathbb{E}|h|\mathbf{1}_{\{2^k s \leq |h| < 2^{k+1} s\}} \leq K \sum_{k=0}^{\infty} 2^{k+1} s \frac{(\log\log 2^k s)^\beta}{2^{2k} s^2}$$
$$= 2K \sum_{k=0}^{\infty} \frac{(\log\log 2^k s)^\beta}{2^k s} = \mathcal{O}\left(\frac{(\log\log s)^\beta}{s}\right). \qquad \square$$

LEMMA 6. *If* $\mathbb{E}(h^2 \wedge u) = \mathcal{O}((\log\log u)^\beta)$ *then*

$$\mathbb{E}\frac{|h|^2}{(\mathrm{LL}\,|h|)^{\beta+\varepsilon}} < \infty$$

*for each* $\varepsilon > 0$.

PROOF. For large $n$,

$$\mathbb{E}\frac{|h|^2}{(\log\log|h|)^{\beta+\varepsilon}}\mathbf{1}_{\{2^{2^{2^n}} \leq |h| < 2^{2^{2^{n+1}}}\}} \leq K \frac{\mathbb{E}(|h|^2 \wedge 2^{2^{2^{n+1}}})}{2^{n(\beta+\varepsilon)}} \leq \tilde{K} \frac{2^{(n+1)\beta}}{2^{n(\beta+\varepsilon)}}$$
$$= \tilde{K} 2^\beta 2^{-n\varepsilon}. \qquad \square$$

**5. The equivalence of several LIL statements.** We will now state general results on the correspondence of the LIL for various kinds of $U$-statistics (as defined in Section 2) based on the same kernel, that we will use extensively in the sequel. Let us start with the following lemma, proved in [5].

LEMMA 7. *Let* $h : \Sigma^d \to \mathbb{R}$ *be a symmetric function. There exist constants* $L_d$, *such that if*

$$(1) \qquad \limsup_{n\to\infty} \frac{1}{(n\log\log n)^{d/2}} \left| \sum_{\mathrm{i}\in I_n^d} h(\mathbf{X}_\mathrm{i}) \right| \leq C \qquad a.s.,$$

*then*

$$(2) \qquad \sum_{n=1}^{\infty} \mathbb{P}\left( \left| \sum_{|\mathrm{i}|\leq 2^n} \epsilon_\mathrm{i}^{\mathrm{dec}} h(\mathbf{X}_\mathrm{i}^{\mathrm{dec}}) \right| > D 2^{nd/2} \log^{d/2} n \right) < \infty$$

*for* $D = L_d C$. *Moreover* (2) *implies*

$$(3) \qquad \limsup_{u\to\infty} \frac{\mathbb{E}(h^2(X) \wedge u)}{(\log\log u)^{d-1}} \leq L_d D^2.$$

LEMMA 8. *For a symmetric function $h:\Sigma^d \to \mathbb{R}$, the LIL (1) is equivalent to the decoupled LIL*

$$\limsup_{n\to\infty} \frac{1}{(n\log\log n)^{d/2}} \left| \sum_{\mathrm{i}\in I_n^d} h(\mathbf{X}_{\mathrm{i}}^{\mathrm{dec}}) \right| \le D \qquad a.s., \tag{4}$$

*meaning that (1) implies (4) with $D = L_d C$, and conversely (4) implies (1) with $C = L_d D$.*

PROOF. We can equivalently write (1) as

$$\lim_{k\to\infty} \mathbb{P}\left( \sup_{n\ge k} \frac{1}{(n\log\log n)^{d/2}} \left| \sum_{\mathrm{i}\in I_n^d} h(\mathbf{X}_{\mathrm{i}}) \right| \ge C+\varepsilon \right) = 0,$$

for all $\varepsilon > 0$, which can be rewritten as

$$\lim_{k\to\infty} \mathbb{P}\left( \left\| \sum_{\substack{|\mathrm{i}|<\infty \\ l\ne j \Rightarrow i_l\ne i_j}} h_{|\mathrm{i}|,k}(\mathbf{X}_{\mathrm{i}}) \right\|_\infty \ge C+\varepsilon \right) = 0, \tag{5}$$

where for $i,k\in\mathbb{N}$, $h_{i,k}$ is an $l^\infty$-valued kernel defined as

$$h_{i,k} = \left( \frac{h}{(k\log\log k)^{d/2}}, \frac{h}{((k+1)\log\log(k+1))^{d/2}}, \dots, \frac{h}{(n\log\log n)^{d/2}}, \dots \right)$$

for $i \le k$ and

$$h_{i,k} = \left( \underbrace{0,\dots,0}_{i-k}, \frac{h}{(i\log\log i)^{d/2}}, \frac{h}{((i+1)\log\log(i+1))^{d/2}}, \dots, \frac{h}{(n\log\log n)^{d/2}}, \dots \right)$$

otherwise. Now the decoupling inequalities by de la Peña and Montgomery–Smith (see [4]) show that (5) is up to constant equivalent to its decoupled version, which is equivalent to (4). □

LEMMA 9. *There exists a universal constant $L<\infty$, such that for any kernel $h:\Sigma_1\times\cdots\times\Sigma_d\to\mathbb{R}$ and variables $(X_i^{(j)})_{i,j}$ like in Definition 3, we have*

$$\mathbb{P}\left( \max_{|\mathrm{j}|\le n} \left| \sum_{\mathrm{i}\,:\, i_k\le j_k, k=1,\dots,d} h(\mathbf{X}_{\mathrm{i}}^{\mathrm{dec}}) \right| \ge t \right) \le L^d \mathbb{P}\left( \left| \sum_{|\mathrm{i}|\le n} h(\mathbf{X}_{\mathrm{i}}^{\mathrm{dec}}) \right| \ge t/L^d \right).$$

PROOF. We will prove by induction with respect to $d$ a stronger statement, namely the inequality in question for Banach space valued $U$-statistics, with the absolute value replaced by the norm. For $d=1$, it is a result by Montgomery–Smith [17]. Assume therefore that the statement holds for kernels of degree smaller than $d$ and consider a kernel $h:\Sigma^d \to B$, for some Banach space $B$. Then, conditioning on $X^{(d)}$, applying the induction assumption to $l_\infty^n(B)$ and $g(x_1,\ldots,x_{d-1}) = (\sum_{i_d \le l} h(x_1,\ldots,x_{d-1},X_{i_d}^{(d)}) : l \le n)$ and finally using the Fubini theorem, we obtain

$$\begin{aligned}
&\mathbb{P}\left(\max_{|\mathrm{j}|\le n}\left\|\sum_{\mathrm{i}\,:\,i_k\le j_k, k=1,\ldots,d} h(\mathbf{X}_{\mathrm{i}}^{\mathrm{dec}})\right\|_B \ge t\right)\\
&\qquad \le L^{d-1}\mathbb{P}\left(\max_{j\le n}\left\|\sum_{|\mathrm{i}|\le n\,:\,i_d\le j} h(\mathbf{X}_{\mathrm{i}}^{\mathrm{dec}})\right\|_B \ge t/L^{d-1}\right).
\end{aligned}$$

Now it is enough to apply the result by Montgomery–Smith, conditionally on $(X^{(1)},\ldots,X^{(d-1)})$. □

COROLLARY 1. *Consider a kernel* $h:\Sigma_1\times\cdots\times\Sigma_d\to\mathbb{R}$, *an array of variables* $(X_i^{(j)})_{i,j}$ *like in Definition 3 and* $\alpha>0$. *If*

$$\sum_{n=1}^{\infty}\mathbb{P}\left(\left|\sum_{|\mathrm{i}|\le 2^n} h(\mathbf{X}_{\mathrm{i}}^{\mathrm{dec}})\right| \ge C2^{n\alpha}\log^\alpha n\right) < \infty,$$

*then*

$$\limsup_{n\to\infty}\frac{1}{(n\log\log n)^\alpha}\left|\sum_{|\mathrm{i}|\le n} h(\mathbf{X}_{\mathrm{i}}^{\mathrm{dec}})\right| \le L_{d,\alpha}C \qquad a.s.$$

PROOF. We have for $0<D<\infty$

$$\begin{aligned}
&\mathbb{P}\left(\sup_{n\ge N}\frac{1}{(n\log\log n)^\alpha}\left|\sum_{|\mathrm{i}|\le n} h(\mathbf{X}_{\mathrm{i}}^{\mathrm{dec}})\right| > D\right)\\
&\qquad \le \mathbb{P}\left(\sup_{k>\lfloor\log N/\log 2\rfloor}\max_{2^{k-1}\le n\le 2^k}\frac{L_\alpha}{(2^k\log k)^\alpha}\left|\sum_{|\mathrm{i}|\le n} h(\mathbf{X}_{\mathrm{i}}^{\mathrm{dec}})\right| > D\right)\\
&\qquad \le \sum_{k>\lfloor\log N/\log 2\rfloor}\mathbb{P}\left(\max_{2^{k-1}\le n\le 2^k}\frac{L_\alpha}{(2^k\log k)^\alpha}\left|\sum_{|\mathrm{i}|\le n} h(\mathbf{X}_{\mathrm{i}}^{\mathrm{dec}})\right| > D\right),
\end{aligned}$$

so the result follows from Lemma 9. □

To prove further statements concerning the equivalence of various types of the LIL, we will have to show that the contribution to a decoupled $U$-statistics from the "diagonal," that is from the sum over multiindices $\mathrm{i}\notin I_n^d$ is negligible. One of our tools will be the following.

LEMMA 10. *If $h : \Sigma^d \to \mathbb{R}$ is canonical and satisfies*

$$\mathbb{E}(h^2 \wedge u) = \mathcal{O}((\log\log u)^\beta),$$

*for some $\beta$, then*

$$\limsup_{n\to\infty} \frac{1}{(n\log\log n)^{d/2}} \Bigg| \sum_{\substack{|\mathrm{i}|\le n \\ \exists_{j\neq k} i_j = i_k}} h(\mathbf{X}_{\mathrm{i}}^{\mathrm{dec}}) \Bigg| = 0 \qquad a.s. \tag{6}$$

PROOF. We will decompose the diagonal into several sums, depending on the "level sets" of the multiindex i. For $\mathcal{J} \in \mathcal{P}_{I_d}$ let $A_{\mathcal{J}}(n)$ be the set of all $|\mathrm{i}| \le n$ such that the index i is constant on all $J \in \mathcal{J}$. Let us notice that the contribution to the sum in (6) from $\mathrm{i} \in A_{\mathcal{J}}(n)$ that is

$$U_{\mathcal{J}}(n) := \sum_{\mathrm{i}\in A_{\mathcal{J}}(n)} h(\mathbf{X}_{\mathrm{i}}^{\mathrm{dec}}),$$

can be treated as a canonical decoupled $U$-statistic of order $\deg \mathcal{J}$ if we only treat the variables $\mathbf{X}_{\mathrm{i}_J}^{\mathrm{dec}}$ as one variable for any $J \in \mathcal{J}$.

Let us now denote for $j < k$, $j, k \in I_d$, $A_{jk} = \{\mathrm{i} : |\mathrm{i}| \le n, i_j = i_k\}$ and $\Delta = \{(j,k) \subseteq I_d^2 : j < k\}$. From the inclusion–exclusion formula we get for every $|\mathrm{i}| \le n$,

$$\mathbf{1}_{\{\exists_{j\neq k} i_j = i_k\}} = \mathbf{1}_{\bigcup_{(j,k)\in\Delta} A_{jk}} = \sum_{l=1}^{\binom{d}{2}} \sum_{\substack{(j_1,k_1),\ldots,(j_l,k_l)\in\Delta \\ \forall_{r\neq s}(j_r,k_r)\neq(j_s,k_s)}} (-1)^{l-1} \mathbf{1}_{A_{j_1k_1}\cap\cdots\cap A_{j_lk_l}}.$$

Hence we have

$$\sum_{\substack{|\mathrm{i}|\le n \\ \exists_{j\neq k} i_j = i_k}} h(\mathbf{X}_{\mathrm{i}}^{\mathrm{dec}}) = \sum_{\substack{\mathcal{J}\in\mathcal{P}_{I_d} \\ \deg\mathcal{J} < d}} a_{\mathcal{J}} U_{\mathcal{J}}(n),$$

for some numbers $a_{\mathcal{J}}$, whose absolute values are bounded by a constant, depending only on $d$. Since the number of summands on the right-hand side does not depend on $n$ either, it is enough to prove that

$$\limsup_{n\to\infty} \frac{|U_{\mathcal{J}}(n)|}{(n\log\log n)^{d/2}} = 0$$

for all $\mathcal{J}$ such that $\deg\mathcal{J} < d$.

Therefore, by Corollary 1, it is enough to prove that for $\deg\mathcal{J} < d$,

$$\sum_{n=1}^{\infty} \mathbb{P}\Bigg( \Bigg| \sum_{\mathrm{i}\in A_{\mathcal{J}}(2^n)} \pi_{\deg\mathcal{J}} h(\mathbf{X}_{\mathrm{i}}^{\mathrm{dec}}) \Bigg| \ge C 2^{nd/2} \log^{d/2} n \Bigg) < \infty \tag{7}$$

for any $C > 0$. (Here $\pi_{\deg \mathcal{J}}$ denotes the Hoeffding projection of the kernel $h$ considered a $U$-statistics of order $\deg \mathcal{J}$, as mentioned above. We have thus actually $\pi_{\deg \mathcal{J}} h = h$.) It is relatively easy to prove (7), as the number of summands is of much smaller order than $2^{nd}$. Obviously $\#A_{\mathcal{J}}(2^n) = 2^{n \deg \mathcal{J}} \le 2^{n(d-1)}$. Let $I$ be any subset of $I_d$, such that for any $J \in \mathcal{J}$, $\#(I \cap J) = 1$. For $h_n = h\mathbf{1}_{\{|h| > 2^{nd}\}}$ we have by Lemma 5

$$\mathbb{E}\left|\sum_{\mathrm{i} \in A_{\mathcal{J}}(2^n)} \epsilon_{\mathrm{i}_I}^{\mathrm{dec}} h_n(\mathbf{X}_{\mathrm{i}}^{\mathrm{dec}})\right| \le 2^{n(d-1)} \mathbb{E}|h_n| \le K 2^{n(d-1)} \frac{\log^\beta n}{2^{nd}} = K \frac{\log^\beta n}{2^n},$$

and the convergence of (7) with $h$ replaced by $h_n$ follows easily from Lemma 1 and the Chebyshev inequality. On the other hand, for $\tilde{h}_n = h\mathbf{1}_{\{|h| \le 2^{nd}\}}$ we have

$$\begin{aligned}\frac{\mathbb{E}|\sum_{\mathrm{i} \in A_{\mathcal{J}}(2^n)} \epsilon_{\mathrm{i}_I}^{\mathrm{dec}} \tilde{h}_n(\mathbf{X}_{\mathrm{i}}^{\mathrm{dec}})|^2}{C^2 2^{nd} \log^d n} &\le \frac{\#A_{\mathcal{J}}(2^n) \mathbb{E}\tilde{h}_n^2}{C^2 2^{nd} \log^d n} \le \frac{2^{n(d-1)} \mathbb{E}\tilde{h}_n^2}{C^2 2^{nd} \log^d n} \\ &\le K C^{-2} 2^{-n} \log^{\beta - d} n,\end{aligned}$$

which (again via Lemma 1 and the Chebyshev inequality) allows us to finish the proof. $\square$

Corollary 2. *The randomized decoupled LIL*

$$\limsup_{n\to\infty} \frac{1}{(n \log\log n)^{d/2}} \left|\sum_{|\mathrm{i}| \le n} \epsilon_{\mathrm{i}}^{\mathrm{dec}} h(\mathbf{X}_{\mathrm{i}}^{\mathrm{dec}})\right| \le C \tag{8}$$

*is equivalent to* (2), *meaning that if* (8) *holds then so does* (2) *with* $D = L_d C$ *and* (2) *implies* (8) *with* $C = L_d D$.

Proof. Implication (2) $\Rightarrow$ (8) follows from Corollary 1. To get (2) from (8), it is enough to show that $\mathbb{E}(h^2 \wedge u) = \mathcal{O}((\log\log u)^d)$, since then by Lemma 10 we can skip the diagonal and by Lemma 8 undecouple to obtain

$$\limsup_{n\to\infty} \frac{1}{(n \log\log n)^{d/2}} \left|\sum_{\mathrm{i} \in I_n^d} \epsilon_{\mathrm{i}} h(\mathbf{X}_{\mathrm{i}})\right| < \infty,$$

which gives (2) by Lemma 7 [note that if $(\varepsilon_i)_i, (\varepsilon_i^{(j)})_i$, $j = 1, \ldots, d$, are independent Rademacher sequences, then so are $(\varepsilon_i \varepsilon_i^{(j)})_i$]. This is, however, easy by a simple modification of arguments from [7], which we will present here for the sake of completeness. Notice that by the Paley–Zygmund inequality and hypercontractivity of Rademacher chaos, we have

$$\mathbb{P}_\varepsilon\left(\left|\sum_{|\mathrm{i}| \le n} \epsilon_{\mathrm{i}}^{\mathrm{dec}} h(\mathbf{X}_{\mathrm{i}}^{\mathrm{dec}})\right| \ge L_d^{-1} \left(\sum_{|\mathrm{i}| \le n} h(\mathbf{X}_{\mathrm{i}}^{\mathrm{dec}})^2\right)^{1/2}\right) \ge \frac{1}{L_d}. \tag{9}$$

Moreover if $\mathbb{E}(h^2 \wedge n) \geq 1$, then

$$\begin{aligned}
&\mathbb{E}\left(\sum_{|\mathrm{i}|\leq n}(h(\mathbf{X}_{\mathrm{i}}^{\mathrm{dec}})^2 \wedge n)\right)^2 \\
&\quad = \sum_{I\subseteq I_d}\sum_{|\mathrm{i}|\leq n}\sum_{\substack{|\mathrm{j}|\leq n \\ \{k\,:\,i_k=j_k\}=I}} \mathbb{E}[h(\mathbf{X}_{\mathrm{i}}^{\mathrm{dec}})^2 \wedge n][h(\mathbf{X}_{\mathrm{j}}^{\mathrm{dec}})^2 \wedge n] \\
&\quad \leq n^{2d}[\mathbb{E}(h^2 \wedge n)]^2 + \sum_{I\subseteq I_d, I\neq\varnothing} n^{d+\#I^c} n\mathbb{E}(h^2 \wedge n) \\
&\quad \leq n^{2d}[\mathbb{E}(h^2 \wedge n)]^2 + (2^d-1)n^{2d}\mathbb{E}(h^2 \wedge n) \\
&\quad \leq 2^d n^{2d}[\mathbb{E}(h^2 \wedge n)]^2 = 2^d\left(\mathbb{E}\sum_{|\mathrm{i}|\leq n}(h(\mathbf{X}_{\mathrm{i}}^{\mathrm{dec}})^2 \wedge n)\right)^2.
\end{aligned}$$

Thus again by Paley–Zygmund, we have

$$\mathbb{P}\left(\sum_{|\mathrm{i}|\leq n} h(\mathbf{X}_{\mathrm{i}}^{\mathrm{dec}})^2 \geq \frac{1}{2}n^d\mathbb{E}(h^2 \wedge n)\right) \geq \frac{1}{L_d},$$

which together with (9) yields

$$\mathbb{P}\left(\left|\sum_{|\mathrm{i}|\leq n} \epsilon_{\mathrm{i}}^{\mathrm{dec}} h(\mathbf{X}_{\mathrm{i}}^{\mathrm{dec}})\right| \geq \tilde{L}_d^{-1} n^{d/2}\sqrt{\mathbb{E}(h^2 \wedge n)}\right) \geq \frac{1}{\tilde{L}_d},$$

which gives $\mathbb{E}(h^2 \wedge n) = \mathcal{O}((\log\log n)^d)$, since by assumption the sequence

$$\frac{1}{(n\log\log n)^{d/2}}\left|\sum_{|\mathrm{i}|\leq n} \epsilon_{\mathrm{i}}^{\mathrm{dec}} h(\mathbf{X}_{\mathrm{i}}^{\mathrm{dec}})\right|$$

is stochastically bounded. □

Corollary 3. *For a symmetric, canonical kernel $h\colon \Sigma^d \to \mathbb{R}$, the LIL (1) is equivalent to the decoupled LIL "with diagonal"*

$$\limsup_{n\to\infty} \frac{1}{(n\log\log n)^{d/2}}\left|\sum_{|\mathrm{i}|\leq n} h(\mathbf{X}_{\mathrm{i}}^{\mathrm{dec}})\right| \leq D \qquad a.s. \tag{10}$$

*again meaning that there are constants $L_d$ such that if (1) holds for some $D$ then so does (10) for $D = L_d C$, and conversely, (10) implies (1) for $C = L_d D$.*

Proof. To show that (1) implies (10) it is enough to use Lemma 8 and then Lemma 10 to add the diagonal (the integrability condition on $h$ follows from Lemma 7).

To obtain the converse implication, it is enough to prove $\mathbb{E}(h(X)^2 \wedge u) = \mathcal{O}((\log\log u)^d)$ since then we are allowed to delete the diagonal by means of Lemma 10 and use Lemma 8 to undecouple the LIL.

From the assumption it follows that for every $\varepsilon > 0$ and sufficiently large $n$,

$$\mathbb{P}\left(\left|\sum_{|\mathrm{i}|\le n} h(\mathbf{X}_{\mathrm{i}}^{\mathrm{dec}})\right| > (D+1)n^{d/2}\log\log^{d/2} n\right) < \varepsilon.$$

Now, by Lemma 9, for arbitrary subsets $A_1,\ldots,A_d \subseteq I_n$,

$$\mathbb{P}\left(\left|\sum_{\mathrm{i}\in A_1\times\cdots\times A_d} h(\mathbf{X}_{\mathrm{i}}^{\mathrm{dec}})\right| > L^d(D+1)n^{d/2}\log\log^{d/2} n\right) \le L^d\varepsilon.$$

Moreover, for fixed values of $(\varepsilon_i^{(j)})$, the expression $\sum_{|\mathrm{i}|\le n}\epsilon_{\mathrm{i}}^{\mathrm{dec}} h(\mathbf{X}_{\mathrm{i}}^{\mathrm{dec}})$ is a sum of $2^d$ expressions of the form $\pm\sum_{\mathrm{i}\in A_1\times\cdots\times A_d} h(\mathbf{X}_{\mathrm{i}}^{\mathrm{dec}})$, where $A_k = \{i : \varepsilon_i^{(k)} = \pm 1\}$. Thus, using the above estimate conditionally, together with the Fubini theorem, we get for sufficiently large $n$,

$$\mathbb{P}\left(\left|\sum_{|\mathrm{i}|\le n}\epsilon_{\mathrm{i}}^{\mathrm{dec}} h(\mathbf{X}_{\mathrm{i}}^{\mathrm{dec}})\right| > 2^d L^d(D+1)n^{d/2}\log\log^{d/2} n\right) \le 2^d L^d\varepsilon.$$

Now we can finish just like in Corollary 2 by applying the Paley–Zygmund inequality and hypercontractive estimates for chaoses. □

**6. The canonical decoupled case.** Before we state the necessary and sufficient conditions for the LIL, let us notice that the integrability condition $\mathbb{E}(h^2 \wedge u) = \mathcal{O}((\log\log u)^{d-1})$ can be equivalently expressed in the language of the $\|\cdot\|_{\mathcal{J},u}$ norms (see Section 3 for the definition). More precisely, we have the following.

LEMMA 11. *For any function $h$ we have*

$$\limsup_{u\to\infty}\frac{(\mathbb{E}(h^2\wedge u))^{1/2}}{(\log\log u)^{(d-1)/2}} = \limsup_{u\to\infty}\frac{\|h\|_{\{I_d\},u}}{(\log\log u)^{(d-1)/2}}.$$

PROOF. Let us denote the lim sup on the right-hand side by $a$, and the other one by $b$. Let us also assume without loss of generality that $h \ge 0$. We will first prove that $a \le b$. Indeed, either $\mathbb{E}(h^2 \wedge u) \le 1$ or we can use

$$f := \frac{h\wedge\sqrt{u}}{(\mathbb{E}(h^2\wedge u))^{1/2}}$$

as a test function in the definition of $\|h\|_{\{I_d\},u}$, thus obtaining for $u \ge 1$

$$\|h\|_{\{I_d\},u} \ge \mathbb{E}hf = \frac{\mathbb{E}(h^2\wedge\sqrt{u}h)}{(\mathbb{E}(h^2\wedge u))^{1/2}} \ge (\mathbb{E}(h^2\wedge u))^{1/2},$$

so we have $(\mathbb{E}(h^2 \wedge u))^{1/2} \le 1 + \|h\|_{\{I_d\},u}$, which immediately yields $a \le b$. To prove the other inequality, let us notice that if $a < \infty$, then for $u$ large enough and any $f$ with $\|f\|_2 \le 1, \|f\|_\infty \le u$ Lemma 5 gives

$$\begin{aligned}\mathbb{E}hf &\le \sqrt{\mathbb{E}h^2\mathbf{1}_{\{h\le u^2\}}} + u\mathbb{E}|h|\mathbf{1}_{\{h\ge u^2\}}\\ &\le (\mathbb{E}(h^2\wedge u^4))^{1/2} + u\frac{K(\log\log u^2)^{d-1}}{u^2},\end{aligned}$$

which gives $b \le a$ since

$$\lim_{u\to\infty}\frac{\log\log u^4}{\log\log u} = 1. \qquad \square$$

THEOREM 5. *Let $h$ be a canonical symmetric kernel in $d$ variables. Then the decoupled LIL*

$$\limsup_{n\to\infty}\frac{1}{n^{d/2}(\log\log n)^{d/2}}\left|\sum_{|\mathrm{i}|\le n} h(\mathbf{X}_{\mathrm{i}}^{\mathrm{dec}})\right| \le C \qquad a.s. \tag{11}$$

*holds if and only if for all $\mathcal{J} \in \mathcal{P}_{I_d}$,*

$$\limsup_{u\to\infty}\frac{1}{(\log\log u)^{(d-\deg\mathcal{J})/2}}\|h\|_{\mathcal{J},u} \le D, \tag{12}$$

*that is, if (11) holds for some $C$ then (12) is satisfied for $D = L_dC$ and conversely, (12) implies (11) with $C = L_dD$.*

PROOF.

*Necessity.* Let us first prove the following.

LEMMA 12. *Let $g:\Sigma^d \to \mathbb{R}$ be a square integrable function. Then*

$$\mathrm{Var}\left(\sum_{|\mathrm{i}|\le n} g(\mathbf{X}_{\mathrm{i}}^{\mathrm{dec}})\right) \le (2^d-1)n^{2d-1}\mathbb{E}g(X)^2.$$

PROOF. We have

$$\begin{aligned}&\mathrm{Var}\left(\sum_{|\mathrm{i}|\le n} g(\mathbf{X}_{\mathrm{i}}^{\mathrm{dec}})\right)\\ &\quad= \mathbb{E}\left(\sum_{|\mathrm{i}|\le n}(g(\mathbf{X}_{\mathrm{i}}^{\mathrm{dec}}) - \mathbb{E}g(\mathbf{X}_{\mathrm{i}}^{\mathrm{dec}}))\right)^2\end{aligned}$$

$$= \sum_{I\subseteq I_d}\sum_{|\mathrm{i}|\le n}\sum_{\substack{|\mathrm{j}|\le n:\\ \{k:\, i_k=j_k\}=I}} \mathbb{E}[(g(\mathbf{X}_{\mathrm{i}}^{\mathrm{dec}})-\mathbb{E}g(\mathbf{X}_{\mathrm{i}}^{\mathrm{dec}}))(g(\mathbf{X}_{\mathrm{j}}^{\mathrm{dec}})-\mathbb{E}g(\mathbf{X}_{\mathrm{j}}^{\mathrm{dec}}))]$$

$$= \sum_{I\subseteq I_d, I\neq\varnothing}\sum_{|\mathrm{i}|\le n}\sum_{\substack{|\mathrm{j}|\le n:\\ \{k:\, i_k=j_k\}=I}} \mathbb{E}[(g(\mathbf{X}_{\mathrm{i}}^{\mathrm{dec}})-\mathbb{E}g(\mathbf{X}_{\mathrm{i}}^{\mathrm{dec}}))(g(\mathbf{X}_{\mathrm{j}}^{\mathrm{dec}})-\mathbb{E}g(\mathbf{X}_{\mathrm{j}}^{\mathrm{dec}}))]$$

$$\le \sum_{I\subseteq I_d, I\neq\varnothing} n^d n^{d-\#I}\,\mathrm{Var}(g(X)) \le (2^d-1)n^{2d-1}\mathbb{E}g(X)^2. \qquad \square$$

Moving to the proof of (12), let us first note that from Corollary 3 and Lemma 7, the series (2) is convergent and (3) holds. Since $\lim_{n\to\infty}\sum_{k=n}^{2n}\frac{1}{k} = \log 2$, there exists $N_0$, such that for all $N > N_0$, there exists $N \le n \le 2N$, satisfying

$$\mathbb{P}\left(\left|\sum_{|\mathrm{i}|\le 2^n}\epsilon_{\mathrm{i}}^{\mathrm{dec}}h(\mathbf{X}_{\mathrm{i}}^{\mathrm{dec}})\right| > L_d C 2^{nd/2}\log^{d/2} n\right) < \frac{1}{10n}. \tag{13}$$

Let us thus fix $N > N_0$ and consider $n$ as above. Let $\mathcal{J} = \{J_1,\ldots,J_k\} \in \mathcal{P}_{I_d}$. Let us also fix functions $f_j\colon \Sigma^{\#J_j}\to\mathbb{R}$, $j=1,\ldots,k$, such that

$$\|f_j(X_{J_j})\|_2 \le 1,$$

$$\|f_j(X_{J_j})\|_\infty \le 2^{n/(2k+1)}.$$

The Chebyshev inequality gives

$$\mathbb{P}\left(\sum_{|\mathrm{i}_{J_j}|\le 2^n} f_j(\mathbf{X}_{\mathrm{i}_{J_j}}^{\mathrm{dec}})^2\log n \le 10\cdot 2^d 2^{\#J_j n}\log n\right) \ge 1 - \frac{1}{10\cdot 2^d}. \tag{14}$$

Moreover, for sufficiently large $N$,

$$\sum_{|\mathrm{i}_{\diamond J_j}|\le 2^n}\frac{1}{2^{n\#J_j}} f_j(\mathbf{X}_{\mathrm{i}_{J_j}}^{\mathrm{dec}})^2\cdot\log n \le \frac{2^{n\#\diamond J_j}2^{2n/(2k+1)}\log n}{2^{n\#J_j}}$$

$$\le \frac{2^{2n/(2k+1)}\log n}{2^n} \le 1.$$

Without loss of generality we may assume that the sequences $(X_i^{(j)})_{i,j}$ and $(\varepsilon_i^{(j)})_{i,j}$ are defined as coordinates of a product probability space. If for each $j=1,\ldots,k$ we denote the set from (14) by $A_k$, we have $\mathbb{P}(\bigcap_{j=1}^k A_k) \ge 0.9$. Recall now Lemma 3. On $\bigcap_{j=1}^k A_k$ we can estimate the $\|\cdot\|_{\mathcal{J},\log n}^*$ norms of the matrix $(h(\mathbf{X}_{\mathrm{i}}^{\mathrm{dec}}))_{|\mathrm{i}|\le 2^n}$ by using the test sequences

$$\alpha_{\mathrm{i}_{J_j}} = \frac{f_j(\mathbf{X}_{\mathrm{i}_{J_j}}^{\mathrm{dec}})\sqrt{\log n}}{10^{1/2}2^{d/2}2^{n\#J_j/2}}.$$

Therefore, with probability at least 0.9, we have

$$\begin{aligned}
&\|(h(\mathbf{X}_{\mathrm{i}}^{\mathrm{dec}}))_{|\mathrm{i}|\le 2^n}\|^*_{\mathcal{J},\log n} \\
&\quad\ge \frac{(\log n)^{k/2}}{2^{dk/2}10^{k/2}2^{(\sum_j \# J_j)n/2}}\Big|\sum_{|\mathrm{i}|\le 2^n} h(\mathbf{X}_{\mathrm{i}}^{\mathrm{dec}})\prod_{j=1}^k f_j(\mathbf{X}_{\mathrm{i}_{J_j}}^{\mathrm{dec}})\Big| \qquad (15)\\
&\quad= \frac{(\log n)^{k/2}}{2^{dk/2}10^{k/2}2^{dn/2}}\Big|\sum_{|\mathrm{i}|\le 2^n} h(\mathbf{X}_{\mathrm{i}}^{\mathrm{dec}})\prod_{j=1}^k f_j(\mathbf{X}_{\mathrm{i}_{J_j}}^{\mathrm{dec}})\Big|.
\end{aligned}$$

Our aim is now to further bound from below the right-hand side of the above inequality, to have, via Lemma 3, control from below on the conditional tail probability of $\sum_{|\mathrm{i}|\le 2^n}\epsilon_{\mathrm{i}}^{\mathrm{dec}}h(\mathbf{X}_{\mathrm{i}}^{\mathrm{dec}})$, given the sample $(X_i^{(j)})$.

From now on let us assume that

$$\Big|\mathbb{E}h(X)\prod_{j=1}^k f_j(X_{J_j})\Big| > 1. \qquad (16)$$

By Corollary 3 and Lemma 7 we have $\mathbb{E}(h^2\wedge u) = \mathcal{O}((\log\log u)^{d-1})$. Thus, the Markov inequality and Lemma 5 give

$$\begin{aligned}
&\mathbb{P}\Big(\Big|\sum_{|\mathrm{i}|\le 2^n} h(\mathbf{X}_{\mathrm{i}}^{\mathrm{dec}})\mathbf{1}_{\{|h(\mathbf{X}_{\mathrm{i}}^{\mathrm{dec}})|>2^n\}}\prod_{j=1}^k f_j(\mathbf{X}_{\mathrm{i}_{J_j}}^{\mathrm{dec}})\Big| \ge \frac{2^{nd}|\mathbb{E}h\prod_{j=1}^k f_j|}{4}\Big)\\
&\quad\le 4\frac{2^{nd}(\prod_{j=1}^k\|f_j\|_\infty)\cdot\mathbb{E}|h|\mathbf{1}_{\{|h|>2^n\}}}{2^{nd}|\mathbb{E}h\prod_{j=1}^k f_j|} \qquad (17)\\
&\quad\le 4\frac{2^{nk/(2k+1)}\mathbb{E}|h|\mathbf{1}_{\{|h|>2^n\}}}{|\mathbb{E}h\prod_{j=1}^k f_j|}\\
&\quad\le 4K\frac{(\log n)^{d-1}}{2^{n(k+1)/(2k+1)}}.
\end{aligned}$$

Let now $h_n = h\mathbf{1}_{\{|h|\le 2^n\}}$. By the Chebyshev inequality, Lemma 12 and (3),

$$\begin{aligned}
&\mathbb{P}\Big(\Big|\sum_{|\mathrm{i}|\le 2^n} h_n(\mathbf{X}_{\mathrm{i}}^{\mathrm{dec}})\prod_{j=1}^k f_j(\mathbf{X}_{\mathrm{i}_{J_j}}^{\mathrm{dec}}) - 2^{nd}\mathbb{E}h_n\prod_{j=1}^k f_j\Big| \ge \frac{2^{nd}}{5}\Big|\mathbb{E}h_n\prod_{j=1}^k f_j\Big|\Big)\\
&\quad\le 25\frac{\mathrm{Var}(\sum_{|\mathrm{i}|\le 2^n} h_n(\mathbf{X}_{\mathrm{i}}^{\mathrm{dec}})\prod_{j=1}^k f_j(\mathbf{X}_{\mathrm{i}_{J_j}}^{\mathrm{dec}}))}{2^{2nd}|\mathbb{E}h_n\prod_{j=1}^k f_j|^2}\\
&\quad\le 25\frac{(2^d-1)2^{n(2d-1)}}{2^{2nd}|\mathbb{E}h_n\prod_{j=1}^k f_j|^2}\mathbb{E}\Big|h_n\prod_{j=1}^k f_j\Big|^2 \qquad (18)
\end{aligned}$$

$$\leq 25(2^d-1)\frac{2^{2nk/(2k+1)}\mathbb{E}h_n^2}{2^n|\mathbb{E}h_n\prod_{j=1}^k f_j|^2}$$

$$\leq 25K(2^d-1)\frac{\log^{d-1} n}{2^{n/(2k+1)}|\mathbb{E}h_n\prod_{j=1}^k f_j|^2}.$$

Let us also notice that for large $n$, by (3), Lemma 5 and (16),

$$\begin{aligned}\Big|\mathbb{E}h_n\prod_{j=1}^k f_j\Big| &\geq \Big|\mathbb{E}h\prod_{j=1}^k f_j\Big| - \Big|\mathbb{E}h\mathbf{1}_{\{|h|>2^n\}}\prod_{j=1}^k f_j\Big| \\ &\geq \Big|\mathbb{E}h\prod_{j=1}^k f_j\Big| - 2^{nk/(2k+1)}K\frac{(\log n)^{d-1}}{2^n} \qquad (19)\\ &\geq \frac{5}{8}\Big|\mathbb{E}h\prod_{j=1}^k f_j\Big| \geq \frac{5}{8}.\end{aligned}$$

Inequalities (17), (18) and (19) imply, that for large $n$ with probability at least 0.9 we have

$$\begin{aligned}&\Big|\sum_{|\mathrm{i}|\leq 2^n} h(\mathbf{X}_{\mathrm{i}}^{\mathrm{dec}})\prod_{j=1}^k f_j(\mathbf{X}_{\mathrm{i}_{J_j}}^{\mathrm{dec}})\Big| \\ &\quad\geq \Big|\sum_{|\mathrm{i}|\leq 2^n} h_n(\mathbf{X}_{\mathrm{i}}^{\mathrm{dec}})\prod_{j=1}^k f_j(\mathbf{X}_{\mathrm{i}_{J_j}}^{\mathrm{dec}})\Big| \\ &\qquad - \Big|\sum_{|\mathrm{i}|\leq 2^n} h(\mathbf{X}_{\mathrm{i}}^{\mathrm{dec}})\mathbf{1}_{\{|h(\mathbf{X}_{\mathrm{i}}^{\mathrm{dec}})|>2^n\}}\prod_{j=1}^k f_j(\mathbf{X}_{\mathrm{i}_{J_j}}^{\mathrm{dec}})\Big| \\ &\quad\geq 2^{nd}\Big(\frac{4}{5}\Big|\mathbb{E}h_n\prod_{j=1}^n f_j\Big| - \frac{1}{4}\Big|\mathbb{E}h\prod_{j=1}^k f_j\Big|\Big) \\ &\quad\geq 2^{nd}\Big(\frac{4}{5}\cdot\frac{5}{8}\Big|\mathbb{E}h\prod_{j=1}^n f_j\Big| - \frac{1}{4}\Big|\mathbb{E}h\prod_{j=1}^k f_j\Big|\Big) \geq \frac{2^{nd}}{4}\Big|\mathbb{E}h\prod_{j=1}^k f_j\Big|.\end{aligned}$$

Together with (15) this yields that for large $n$ with probability at least 0.8,

$$\|(h_{\mathrm{i}})_{|\mathrm{i}|\leq 2^n}\|^*_{\mathcal{J},\log n} \geq \frac{2^{nd/2}\log^{k/2} n}{4\cdot 2^{dk/2}10^{k/2}}\Big|\mathbb{E}h\prod_{j=1}^k f_j\Big|.$$

Thus, by Lemma 3, for large $n$

$$\mathbb{P}\left(\left|\sum_{|\mathrm{i}|\leq 2^n}\epsilon_{\mathrm{i}}^{\mathrm{dec}}h(\mathbf{X}_{\mathrm{i}}^{\mathrm{dec}})\right|\geq c_d\frac{2^{nd/2}\log^{k/2}n}{4\cdot 2^{dk/2}10^{k/2}}\left|\mathbb{E}h\prod_{j=1}^k f_j\right|\right)\geq\frac{8}{10n},$$

which together with (13) gives

$$\left|\mathbb{E}h\prod_{j=1}^k f_j\right|\leq L_d C\frac{4\cdot 2^{dk/2}10^{k/2}}{c_d}\log^{(d-k)/2}n.$$

In particular for sufficiently large $N$, for arbitrary functions $f_j:\Sigma^{\#J_j}\to\mathbb{R}$, $j=1,\dots,k$, such that

$$\|f_j(X_{J_j})\|_2\leq 1,\qquad \|f_j(X_{J_j})\|_\infty\leq 2^{N/(2k+1)}$$

we have

$$\left|\mathbb{E}h\prod_{j=1}^k f_j\right|\leq L_d C\frac{4\cdot 2^{dk/2}10^{k/2}}{c_d}\log^{(d-k)/2}n\leq\tilde{L}_d C\log^{(d-k)/2}N.$$

Thus, for large $u$ ($u\geq u_0$),

$$\sup\left\{\left|\mathbb{E}h(X)\prod_{j=1}^k f_j(X_{J_j})\right|:\|f_j(X_{J_j})\|_2\leq 1,\|f_j(X_{J_j})\|_\infty\leq u^{1/(2k+1)}\right\}$$
$$\leq\bar{L}_d(\log\log u)^{(d-k)/2},$$

and so

$$\sup\left\{\left|\mathbb{E}h(X)\prod_{j=1}^k f_j(X_{J_j})\right|:\|f_j(X_{J_j})\|_2\leq 1,\|f_j(X_{J_j})\|_\infty\leq u\right\}$$
$$\leq\hat{L}_d(\log\log u)^{(d-k)/2},$$

for all $u\geq u_0^{1/(2k+1)}$, which proves the necessity part of the theorem.

*Sufficiency.* The proof consists of several truncation arguments. In the first part, until the $\|\cdot\|_{\mathcal{J},u}$ norms come into play, we follow the lines of the proof of the special case $d=2$, presented in [8], with some modifications. At each step we will show that

$$\sum_{n=1}^\infty\mathbb{P}\left(\left|\sum_{|\mathrm{i}|\leq 2^n}\pi_d h_n(\mathbf{X}_{\mathrm{i}}^{\mathrm{dec}})\right|\geq C2^{nd/2}\log^{d/2}n\right)<\infty,\tag{20}$$

with $h_n=h\mathbf{1}_{A_n}$ for some sequence of sets $A_n$.

*Step* 1. Inequality (20) holds for any $C > 0$ if

$$A_n \subseteq \{x : h^2(x) \geq 2^{nd} \log^d n\}.$$

We have, by the Chebyshev inequality and the inequality $\mathbb{E}|\pi_d h_n| \leq 2^d \mathbb{E}|h_n|$ (which follows directly from the definition of $\pi_d$ or may be considered a trivial case of Lemma 1),

$$\begin{aligned}
&\sum_n \mathbb{P}\left(\left|\sum_{|\mathrm{i}| \leq 2^n} \pi_d h_n(\mathbf{X}_{\mathrm{i}}^{\mathrm{dec}})\right| \geq C 2^{nd/2} \log^{d/2} n\right) \\
&\quad \leq \sum_n \frac{\mathbb{E}|\sum_{|\mathrm{i}| \leq 2^n} \pi_d h_n(\mathbf{X}_{\mathrm{i}}^{\mathrm{dec}})|}{C 2^{nd/2} \log^{d/2} n} \\
&\quad \leq 2^d \sum_n \frac{2^{nd} \mathbb{E}|h| \mathbf{1}_{\{|h| \geq 2^{nd/2} \log^{d/2} n\}}}{C 2^{nd/2} \log^{d/2} n} \\
&\quad = 2^d C^{-1} \mathbb{E}|h| \sum_n \frac{2^{nd/2}}{\log^{d/2} n} \mathbf{1}_{\{|h| \geq 2^{nd/2} \log^{d/2} n\}} \\
&\quad \leq L_d C^{-1} \mathbb{E} \frac{|h|^2}{(\mathrm{LL}\,|h|)^d} < \infty,
\end{aligned}$$

where the last inequality follows from Lemma 6, Lemma 11 and condition (12) for $\mathcal{J} = \{I_d\}$.

*Step* 2. Inequality (20) holds for any $C > 0$ if

$$A_n \subseteq \{x \in \Sigma^d : h^2(x) \leq 2^{2nd}, \exists_{I \neq \varnothing, I_d} \mathbb{E}_I(h^2 \wedge 2^{2nd}) \geq 2^{\# I^c n} \log^d n\}.$$

By Lemma 1 and the Chebyshev inequality, it is enough to prove that

$$\sum_n \frac{\mathbb{E}|\sum_{|\mathrm{i}| \leq 2^n} \epsilon_{\mathrm{i}}^{\mathrm{dec}} h_n(\mathbf{X}_{\mathrm{i}}^{\mathrm{dec}})|}{2^{nd/2} \log^{d/2} n} < \infty.$$

The set $A_n$ can be written as

$$\bigcup_{I \subseteq I_d, I \neq I_d, \varnothing} A_n(I),$$

where the sets $A_n(I)$ are pairwise disjoint and

$$A_n(I) \subseteq \{x : h^2(x) \leq 2^{2nd}, \mathbb{E}_I(h^2 \wedge 2^{2nd}) \geq 2^{\# I^c n} \log^d n\}.$$

Therefore, it suffices to prove that

$$\sum_n \frac{\mathbb{E}|\sum_{|\mathrm{i}| \leq 2^n} \epsilon_{\mathrm{i}}^{\mathrm{dec}} h(\mathbf{X}_{\mathrm{i}}^{\mathrm{dec}}) \mathbf{1}_{A_n(I)}(\mathbf{X}_{\mathrm{i}}^{\mathrm{dec}})|}{2^{nd/2} \log^{d/2} n} < \infty. \tag{21}$$

Let for $l \in \mathbb{N}$,

$$A_{n,l}(I) := \{x : h^2(x) \le 2^{2nd},\\ 2^{2l+2+\#I^c n}\log^d n > \mathbb{E}_I(h^2 \wedge 2^{2nd}) \ge 2^{2l+\#I^c n}\log^d n\} \cap A_n(I).$$

Then $h_n \mathbf{1}_{A_n(I)} = \sum_{l=0}^{\infty} h_{n,l}$, where $h_{n,l} := h_n \mathbf{1}_{A_{n,l}(I)}$.

We have

$$\begin{aligned}
&\mathbb{E}\Big|\sum_{|\mathrm{i}|\le 2^n} \epsilon_{\mathrm{i}}^{\mathrm{dec}} h_{n,l}(\mathbf{X}_{\mathrm{i}}^{\mathrm{dec}})\Big| \\
&\quad \le \sum_{|\mathrm{i}_{I^c}|\le 2^n} \mathbb{E}_{I^c}\mathbb{E}_I\Big|\sum_{|\mathrm{i}_I|\le 2^n} \epsilon_{\mathrm{i}_I}^{\mathrm{dec}} h_{n,l}(\mathbf{X}_{\mathrm{i}}^{\mathrm{dec}})\Big| \\
&\quad \le \sum_{|\mathrm{i}_{I^c}|\le 2^n} \mathbb{E}_{I^c}\bigg(\mathbb{E}_I\Big|\sum_{|\mathrm{i}_I|\le 2^n} \epsilon_{\mathrm{i}_I}^{\mathrm{dec}} h_{n,l}(\mathbf{X}_{\mathrm{i}}^{\mathrm{dec}})\Big|^2\bigg)^{1/2} \\
&\quad \le 2^{(\#I^c+\#I/2)n}\mathbb{E}_{I^c}(\mathbb{E}_I|h_{n,l}|^2)^{1/2} \\
&\quad \le 2^{(\#I^c+d/2)n+l+1}\log^{d/2} n\,\mathbb{P}_{I^c}(\mathbb{E}_I(h^2\wedge 2^{2nd}) \ge 2^{2l+\#I^c n}\log^d n).
\end{aligned}$$

Therefore, to get (21), it is enough to show that

$$\sum_{l=0}^{\infty}\sum_n 2^{l+\#I^c n}\mathbb{P}_{I^c}(\mathbb{E}_I(h^2\wedge 2^{2nd}) \ge 2^{2l+\#I^c n}\log^d n) < \infty.$$

But this is just the statement of Lemma 4 for $a = 2d$.

*Step* 3. Inequality (20) holds for any $C > 0$ if

$$A_n \subseteq \{x : 2^{nd} n^{-2d} < h^2(x) \le 2^{nd}\log^d n\\ \text{and } \forall_{I\neq\varnothing, I_d}\ \mathbb{E}_I(h^2\wedge 2^{2nd}) \le 2^{\#I^c n}\log^d n\}.$$

By Lemma 1 and the Chebyshev inequality, it is enough to show that

$$\sum_n \frac{\mathbb{E}|\sum_{|\mathrm{i}|\le 2^n}\epsilon_{\mathrm{i}}^{\mathrm{dec}} h_n(\mathbf{X}_{\mathrm{i}}^{\mathrm{dec}})|^4}{2^{2nd}\log^{2d} n} < \infty.$$

The Khintchine inequality for Rademacher chaoses gives

$$\begin{aligned}
L_d^{-1}\mathbb{E}\Big|\sum_{|\mathrm{i}|\le 2^n}\epsilon_{\mathrm{i}}^{\mathrm{dec}} h_n(\mathbf{X}_{\mathrm{i}}^{\mathrm{dec}})\Big|^4 &\le \mathbb{E}\bigg(\sum_{|\mathrm{i}|\le 2^n} h_n(\mathbf{X}_{\mathrm{i}}^{\mathrm{dec}})^2\bigg)^2 \\
&= \sum_{I\subseteq I_d}\sum_{|\mathrm{i}|\le 2^n}\sum_{\substack{|\mathrm{j}|\le 2^n:\\ \{k:\, i_k=j_k\}=I}} \mathbb{E}h_n(\mathbf{X}_{\mathrm{i}}^{\mathrm{dec}})^2 h_n(\mathbf{X}_{\mathrm{j}}^{\mathrm{dec}})^2
\end{aligned}$$

$$\leq \sum_{I\subseteq I_d} 2^{nd}2^{n(d-\#I)}\mathbb{E}[h_n(\mathbf{X})^2\cdot h_n(\tilde{\mathbf{X}}(I))^2],$$

where $\mathbf{X}=(X_1,\ldots,X_d)$ and $\tilde{\mathbf{X}}(I)=((X_i)_{i\in I},(X_i^{(1)})_{i\in I^c})$.

To prove the statement of this step it thus suffices to show that for $I\subseteq I_d$ we have

$$\sum_n \frac{2^{-n\#I}}{\log^{2d} n}\mathbb{E}[h_n(X)^2 h_n(\tilde{X}(I))^2]<\infty.$$

(a) $I=I_d$. Then

$$\sum_n \frac{\mathbb{E}h_n^4}{2^{nd}\log^{2d} n}\leq \mathbb{E}h^4\sum_n \frac{1}{2^{nd}\log^{2d} n}\mathbf{1}_{\{h^2\leq 2^{nd}\log^d n\}}$$
$$\leq L_d\mathbb{E}h^4\frac{1}{h^2(\mathrm{LL}\,|h|)^d}<\infty$$

by Lemma 6.

(b) $I\neq I_d,\varnothing$. Let us denote by $\mathbb{E}_I,\mathbb{E}_{I^c},\tilde{\mathbb{E}}_{I^c}$ respectively the expectation with respect to $(X_i)_{i\in I}$, $(X_i)_{i\in I^c}$ and $(X_i^{(1)})_{i\in I^c}$. Let also $\tilde{h}$, $\tilde{h}_n$ stand respectively for $h(\tilde{\mathbf{X}}(I))$, $h_n(\tilde{\mathbf{X}}(I))$. Then

$$\sum_n \frac{\mathbb{E}(h_n^2\cdot\tilde{h}_n^2)}{2^{n\#I}\log^{2d} n}$$
$$\leq 2\sum_n \frac{\mathbb{E}(h_n^2\cdot\tilde{h}_n^2\mathbf{1}_{\{|h|\leq|\tilde{h}|\}})}{2^{n\#I}\log^{2d} n}$$
$$\leq 2\mathbb{E}h^2\tilde{h}^2\mathbf{1}_{\{|h|\leq|\tilde{h}|\}}\sum_n \frac{1}{2^{n\#I}\log^{2d} n}\mathbf{1}_{\{\mathbb{E}_{I^c}(h^2\wedge 2^{2nd})\leq 2^{\#In}\log^d n,\tilde{h}^2\leq 2^{2nd}\}}$$
$$\leq 2\mathbb{E}h^2\tilde{h}^2\mathbf{1}_{\{|h|\leq|\tilde{h}|\}}\sum_n \frac{1}{2^{n\#I}\log^{2d} n}\mathbf{1}_{\{\mathbb{E}_{I^c}(h^2\wedge\tilde{h}^2)\leq 2^{\#In}\log^d n,\tilde{h}^2\leq 2^{2nd}\}}$$
$$\leq L_d\mathbb{E}h^2\tilde{h}^2\mathbf{1}_{\{|h|\leq|\tilde{h}|\}}\frac{1}{(\mathbb{E}_{I^c}(h^2\wedge\tilde{h}^2))(\mathrm{LL}|\tilde{h}|)^d}$$
$$= L_d\mathbb{E}_I\tilde{\mathbb{E}}_{I^c}\tilde{h}^2\mathbb{E}_{I^c}h^2\mathbf{1}_{\{|h|\leq|\tilde{h}|\}}\frac{1}{(\mathbb{E}_{I^c}(h^2\wedge\tilde{h}^2))(\mathrm{LL}|\tilde{h}|)^d}$$
$$\leq L_d\mathbb{E}\frac{\tilde{h}^2}{(\mathrm{LL}|\tilde{h}|)^d}<\infty$$

by Lemma 6.

(c) $I = \varnothing$. We have,

$$\sum_n \frac{(\mathbb{E}h_n^2)^2}{\log^{2d} n} \leq K \sum_n \frac{\mathbb{E}h_n^2}{\log^{d+1} n}$$
(22)
$$\leq K\mathbb{E}h^2 \sum_n \frac{1}{\log^{d+1} n} \mathbf{1}_{\{2^{nd}n^{-2d} < h^2 \leq 2^{nd}\log^d n\}}.$$

For $M > 0$ let us now estimate $\#\{n : 2^{nd}n^{-2d} \leq M \leq 2^{nd}(\log n)^d\}$. Let $n_{\max}, n_{\min}$ denote the greatest and the smallest element of this set. Then

$$n_{\min} \log 2 + \log\log n_{\min} \geq \frac{\log M}{d},$$
$$n_{\max} \log 2 - 2\log n_{\max} \leq \frac{\log M}{d},$$

hence

$$(n_{\max} - n_{\min}) \log 2 \leq 2\log n_{\max} + \log\log n_{\min} \leq 3\log n_{\max}$$
$$\leq L \log\log M.$$

The right-hand side of (22) is thus bounded by

$$K\mathbb{E}\frac{|h|^2 \mathrm{LL}|h|}{(\mathrm{LL}|h|)^{d+1}} = K\mathbb{E}\frac{h^2}{(\mathrm{LL}|h|)^d} < \infty$$

by Lemma 6.

*Step* 4. Inequality (20) holds for any $C > 0$ if

$$A_n \subseteq \Big\{x : h^2 \leq 2^{nd}n^{-2d}, \forall_{I \neq \varnothing, I_d} \mathbb{E}_I(h^2 \wedge 2^{2nd}) \leq 2^{\#I^c n} \log^d n,$$
$$\exists_{I \neq \varnothing, I_d} \frac{2^{\#I^c n}}{\log^{\#I^c} n} \leq \mathbb{E}_I(h^2 \wedge 2^{2nd})\Big\}.$$

The only difference between this step and the previous one is the proof of convergence in the case (c), as in the two other cases we were using only bounds from above on $h^2$ and $\mathbb{E}_I(h^2 \wedge 2^{2nd})$, which are still valid.

Let us notice, that

$$\mathbb{E}h_n^2 \leq \sum_{I \subseteq I_d, I \neq \varnothing, I_d} \mathbb{E}(h^2 \wedge 2^{2nd})\mathbf{1}_{\{(\log n)^{-\#I^c} \leq 2^{-\#I^c n}\mathbb{E}_I(h^2 \wedge 2^{2nd}) \leq (\log n)^d\}}$$
$$\leq \sum_{I \subseteq I_d, I \neq \varnothing, I_d} \sum_{k=1}^{d+\#I^c} \mathbb{E}_{I^c}\mathbb{E}_I(h^2 \wedge 2^{2nd})\mathbf{1}_{\{(\log n)^{d-k} \leq 2^{-\#I^c n}\mathbb{E}_I(h^2 \wedge 2^{2nd}) \leq (\log n)^{d+1-k}\}}$$
$$\leq \sum_{I \subseteq I_d, I \neq \varnothing, I_d} \sum_{k=1}^{d+\#I^c} 2^{\#I^c n}(\log n)^{d+1-k}\mathbb{P}_{I^c}(\mathbb{E}_I(h^2 \wedge 2^{2nd}) \geq 2^{\#I^c n}(\log n)^{d-k}).$$

Thus

$$\sum_n \frac{2^{2nd}(\mathbb{E}h_n^2)^2}{2^{2nd}(\log n)^{2d}}$$
$$\leq \tilde{K}\sum_n \frac{\mathbb{E}h_n^2}{(\log n)^{d+1}}$$
$$\leq K \sum_{I\subseteq I_d, I\neq\varnothing, I_d} \sum_{k=1}^{d+\#I^c} \sum_n \frac{2^{\#I^c n}}{(\log n)^k}\mathbb{P}_{I^c}(\mathbb{E}_I(h^2\wedge 2^{2nd}) \geq 2^{\#I^c n}(\log n)^{d-k})$$
$$< \infty$$

by Lemma 4.

*Step* 5. Inequality (20) holds for *some* $C \leq L_d D$ if

$$A_n = \left\{x : h^2 \leq 2^{nd}n^{-2d}, \forall_{I\neq\varnothing, I_d}\ \mathbb{E}_I(h^2\wedge 2^{2nd}) \leq \frac{2^{\#I^c n}}{\log^{\#I^c} n}\right\}.$$

This is the only part of the proof in which we use the assumptions on the $\|\cdot\|_{\mathcal{J},u}$ norms of $h$ for $\deg\mathcal{J} > 1$. Our aim is to estimate $\|h_n\|_{\mathcal{J}}$ and then use Theorem 4.

Let us note that we can assume that

$$D = 1 \tag{23}$$

[if $D \neq 0$ then we simply scale the function, otherwise (12) for $\mathcal{J} = \{\{1\},\dots,\{d\}\}$ gives $h = 0$].

Let us thus consider $\mathcal{J} = \{J_1,\dots,J_k\} \in \mathcal{P}_{I_d}$ and denote as usual $X = (X_1,\dots,X_d)$, $X_I = (X_i)_{i\in I}$. Recall that

$$\|h_n\|_{\mathcal{J}} = \sup\left\{\mathbb{E}\left[h_n(X)\prod_{i=1}^k f_i(X_{J_i})\right] : \mathbb{E}f_i^2(X_{J_i}) \leq 1\right\}.$$

In what follows, to simplify the already quite complicated notation, let us suppress the arguments of all the functions and write just $h$ instead of $h(X)$ and $f_i$ instead of $f_i(X_{J_i})$.

Let us notice that if $\mathbb{E}f_i^2 \leq 1$, $i = 1,\dots,k$, then for each $j = 1,\dots,k$ and $J \subsetneq J_j$ by the Schwarz inequality applied conditionally to $X_{J_j\setminus J}$,

$$\mathbb{E}\left|h_n\prod_{i=1}^k f_i \mathbf{1}_{\{\mathbb{E}_J f_j^2 > a^2\}}\right| \leq \mathbb{E}_{J_j\setminus J}\left[\left(\mathbb{E}_{(J_j\setminus J)^c}\prod_{i=1}^k f_i^2\right)^{1/2}\mathbf{1}_{\{\mathbb{E}_J f_j^2 \geq a^2\}}(\mathbb{E}_{(J_j\setminus J)^c}h_n^2)^{1/2}\right]$$
$$\leq \mathbb{E}_{J_j\setminus J}[(\mathbb{E}_J f_j^2)^{1/2}\mathbf{1}_{\{\mathbb{E}_J f_j^2 \geq a^2\}}(\mathbb{E}_{(J_j\setminus J)^c}h_n^2)^{1/2}]$$

$$\leq 2^{n\#(J_j\setminus J)/2}\mathbb{E}_{J_j\setminus J}[(\mathbb{E}_J f_j^2)^{1/2}\mathbf{1}_{\{\mathbb{E}_J f_j^2\geq a^2\}}]$$
$$\leq 2^{n\#(J_j\setminus J)/2}a^{-1}.$$

This way we obtain

$$\|h_n\|_{\mathcal{J}} \leq \sup\bigg\{\mathbb{E}\bigg[h_n\prod_{i=1}^k f_i\bigg] : \|f_i\|_2\leq 1,$$
$$\|(\mathbb{E}_J f_i^2)^{1/2}\|_\infty \leq 2^{n\#(J_i\setminus J)/2} \text{ for } J\subsetneq J_i\bigg\}$$
$$+\sum_{i=1}^k(2^{\#J_i}-1) \tag{24}$$
$$\leq L_d+\sup\bigg\{\mathbb{E}\bigg[h_n\prod_{i=1}^k f_i\bigg] : \|f_i\|_2\leq 1,$$
$$\|(\mathbb{E}_J f_i^2)^{1/2}\|_\infty \leq 2^{n\#(J_i\setminus J)/2} \text{ for } J\subsetneq J_i\bigg\}.$$

Let us thus consider arbitrary $f_i$, $i=1,\ldots,k$, such that $\|f_i\|_2\leq 1$, $\|(\mathbb{E}_J f_i^2)^{1/2}\|_\infty \leq 2^{n\#(J_i\setminus J)/2}$ for $J\subsetneq J_i$ (note that the latter condition means in particular that $\|f_i\|_\infty\leq 2^{n\#J_i/2}$).

We have, by assumptions (12) and (23) for large $n$,

$$\bigg|\mathbb{E}\bigg[h\prod_{i=1}^k f_i\bigg]\bigg| \leq \|h\|_{\mathcal{J},2^{nd/2}} \leq L_d\log^{(d-\deg\mathcal{J})/2} n. \tag{25}$$

For sufficiently large $n$,

$$\mathbb{E}\bigg|h\mathbf{1}_{\{|h|\geq 2^{nd/2}n^d\}}\prod_{i=1}^k f_i\bigg| \leq \mathbb{E}|h|\mathbf{1}_{\{|h|\geq 2^{nd/2}n^d\}}\prod_{i=1}^k\|f_i\|_\infty \leq K2^{nd/2}\frac{\log^{d-1}n}{2^{nd/2}n^d}\leq 1,$$

where the second inequality follows from Lemma 5.

Moreover, if we denote $\tilde{h}_n=|h|\wedge 2^{d\cdot\exp(\lceil\log n\rceil)}$, we get for $I\subseteq I_d, I\neq\varnothing, I_d$,

$$\mathbb{E}\bigg|\tilde{h}_n\prod_{i=1}^k f_i\mathbf{1}_{\{\mathbb{E}_I\tilde{h}_n^2\geq 2^{n\#I^c}n\}}\bigg| \leq \mathbb{E}_{I^c}\bigg[(\mathbb{E}_I\tilde{h}_n^2)^{1/2}\mathbf{1}_{\{\mathbb{E}_I\tilde{h}_n^2\geq 2^{n\#I^c}n\}}\prod_{i=1}^k(\mathbb{E}_{J_i\cap I}f_i^2)^{1/2}\bigg]$$
$$\leq\prod_{i=1}^k 2^{n\#(J_i\cap I^c)/2}\mathbb{E}_{I^c}[(\mathbb{E}_I\tilde{h}_n^2)^{1/2}\mathbf{1}_{\{\mathbb{E}_I\tilde{h}_n^2\geq 2^{n\#I^c}n\}}]$$
$$\leq 2^{n\#I^c/2}\frac{\mathbb{E}\tilde{h}_n^2}{2^{n\#I^c/2}\sqrt{n}}\leq K\frac{\log^{d-1}n}{\sqrt{n}}\leq 1$$

for large $n$.

By the last three inequalities we obtain

$$\begin{aligned}
&\left|\mathbb{E}\left[h_n \prod_{i=1}^k f_i\right]\right| \\
&\quad \le \left|\mathbb{E}h \prod_{i=1}^k f_i\right| + \left|\mathbb{E}h\mathbf{1}_{A_n^c} \prod_{i=1}^k f_i\right| \\
&\quad \le L_d \log^{(d-\deg\mathcal{J})/2} n + \mathbb{E}\left|h \prod_{i=1}^k f_i \mathbf{1}_{\{|h|\ge 2^{nd/2}n^{-d}\}}\right| \\
&\qquad + \sum_{I\subseteq I_d, I\ne\varnothing, I_d} \mathbb{E}\left|h\mathbf{1}_{\{|h|<2^{nd/2}n^{-d}\}} \prod_{i=1}^k f_i \mathbf{1}_{\{\mathbb{E}_I(h^2\wedge 2^{2nd})\ge 2^{n\#I^c}(\log n)^{-\#I^c}\}}\right| \\
&\quad \le L_d \log^{(d-\deg\mathcal{J})/2} n + 1 + (\mathbb{E}|h|^2 \mathbf{1}_{\{2^{nd/2}n^{-d}\le|h|\le 2^{nd/2}n^d\}})^{1/2} \\
&\qquad + \sum_{I\subseteq I_d, I\ne\varnothing, I_d} \mathbb{E}\left|\tilde{h}_n \prod_{i=1}^k f_i \mathbf{1}_{\{\mathbb{E}_I \tilde{h}_n^2 \ge 2^{n\#I^c}(\log n)^{-\#I^c}\}}\right| \\
&\quad \le L_d \log^{(d-\deg\mathcal{J})/2} n + 2^d + \sum_{I\subsetneq I_d} d_n(I),
\end{aligned}$$

where

$$\begin{aligned}
d_n(I)^2 &= \mathbb{E}\tilde{h}_n^2 \mathbf{1}_{\{2^{n\#I^c}n^{-1}\le \mathbb{E}_I\tilde{h}_n^2\le 2^{n\#I^c}n\}} \qquad \text{for } I\ne\varnothing, I_d, \\
d_n(\varnothing)^2 &= \mathbb{E}h^2\mathbf{1}_{\{2^{nd/2}n^{-d}\le|h|\le 2^{nd/2}n^d\}}.
\end{aligned}$$

Using (24) we eventually obtain

$$\|h_n\|_{\mathcal{J}} \le L_d \log^{(d-\deg\mathcal{J})/2} n + D_n, \tag{26}$$

where $D_n = \sum_{I\subsetneq I_d} d_n(I)$.

This estimate will allow us to finish the proof by means of the following

Lemma 13. *For sufficiently large $C = L_d A$ and all $\mathcal{J}\in\mathcal{P}_{I_d}$,*

$$\sum_n \exp\left(-\left(\frac{C\log^{d/2} n}{A(\log^{(d-\deg\mathcal{J})/2} n + D_n)}\right)^{2/\deg\mathcal{J}}\right) < \infty.$$

Proof. Let us notice that for $k=1,2,\dots$

$$\sum_{k<\log n\le k+1} h^2\mathbf{1}_{\{2^{nd/2}n^{-d}\le|h|\le 2^{nd/2}n^d\}} \le L_d(k+1)^2\mathbf{1}_{\{|h|\le 2^{ne^{k+1}/2}e^{d(k+1)}\}}$$

and

$$\sum_{k<\log n\le k+1} \mathbb{E}_I \tilde{h}_n^2 \mathbf{1}_{\{2^{n\# I^c} n^{-1}\le \mathbb{E}_I \tilde{h}_n^2 \le 2^{n\# I^c} n\}}$$
$$= \sum_{k<\log n\le k+1} \mathbb{E}_I \tilde{h}_{e^{k+1}}^2 \mathbf{1}_{\{2^{n\# I^c} n^{-1}\le \mathbb{E}_I \tilde{h}_{e^{k+1}}^2 \le 2^{n\# I^c} n\}}$$
$$\le L_d \mathbb{E}_I \tilde{h}_{e^{k+1}}^2 (k+1) = L_d(k+1)\mathbb{E}_I(h^2 \wedge 2^{2de^{k+1}}),$$

since for any numbers $1\le a,b\le d$ and $x\ge 0$, the number of intervals of the form $[2^{na}n^{-b}, 2^{na}n^{b}]$ with $k<\log n\le k+1$, containing $x$ is smaller than $L_d(k+1)$.

Integrating the above inequalities and using Lemma 11, assumption (12) for $\mathcal{J}=\{I_d\}$, assumption (23) and the Cauchy–Schwarz inequality we get

$$\sum_{k<\log n\le k+1} D_n^2 \le (2^d-1) \sum_{k<\log n\le k+1} \sum_{I\subsetneq I_d} d_n(I)^2 \le L_d(k+1)^d.$$

Thus

$$\#\{n : k<\log n\le k+1, D_n\ge 1\} \le L_d(k+1)^d$$

and therefore for $C$ large enough (since $D_n \le L_d \log^{(d-1)/2} n$)

$$\sum_n \exp\Big(-\Big(\frac{C\log^{d/2} n}{A(\log^{(d-\deg\mathcal{J})/2} n + D_n)}\Big)^{2/\deg\mathcal{J}}\Big)$$
$$\le \sum_n \exp(-(C/2A)^{2/\deg\mathcal{J}} \log n)$$
$$+ \sum_{D_n\ge 1} \exp\Big(\Big(-\frac{C\log^{d/2} n}{A(1+L_d)\log^{(d-1)/2} n}\Big)^{2/\deg\mathcal{J}}\Big)$$
$$\le \sum_n \exp(-(C/2A)^{2/\deg\mathcal{J}} \log n)$$
$$+ L_d \sum_k (k+1)^d \exp(-(C/A(1+L_d))^{2/\deg\mathcal{J}} k^{1/\deg\mathcal{J}}) < \infty$$

for $C=\tilde{L}_d A$. □

Going back to the proof of Step 5, let us notice that by Theorem 4 and (26), we have

$$\mathbb{P}\Big(\Big|\sum_{|\mathbf{i}|\le 2^n} \pi_d h_n(\mathbf{X}_{\mathbf{i}}^{\mathrm{dec}})\Big| \ge C 2^{nd/2}\log^{d/2} n\Big)$$

$$\begin{aligned}
&\le L_d \sum_{\mathcal{J}\in\mathcal{P}_{I_d}} \exp\biggl(L_d^{-1}\biggl(\frac{C2^{nd/2}\log^{d/2} n}{2^{nd/2}\|h_n\|_{\mathcal{J}}}\biggr)^{2/\deg\mathcal{J}}\biggr)\\
&\qquad + L_d \sum_{I\subsetneq I_d} \exp\biggl(L_d^{-1}\biggl(\frac{C2^{nd/2}\log^{d/2} n}{2^{n\#I/2}\|(\mathbb{E}_I h_n^2)^{1/2}\|_\infty}\biggr)^{2/(d+\#I^c)}\biggr)\\
&\le L_d \sum_{\mathcal{J}\in\mathcal{P}_{I_d}} \exp\biggl(L_d^{-1}\biggl(\frac{C\log^{d/2} n}{L_d\log^{(d-\deg\mathcal{J})/2} n + D_n}\biggr)^{2/\deg\mathcal{J}}\biggr)\\
&\qquad + L_d \sum_{I\subsetneq I_d} \exp\biggl(L_d^{-1}\biggl(\frac{C2^{nd/2}\log^{d/2} n}{2^{n\#I/2}2^{n\#I^c/2}\log^{-\#I^c/2} n}\biggr)^{2/(d+\#I^c)}\biggr)\\
&\le L_d \sum_{\mathcal{J}\in\mathcal{P}_{I_d}} \exp\biggl(L_d^{-1}\biggl(\frac{C\log^{d/2} n}{L_d\log^{(d-\deg\mathcal{J})/2} n + D_n}\biggr)^{2/\deg\mathcal{J}}\biggr)\\
&\qquad + L_d \sum_{I\subsetneq I_d} \exp(L_d^{-1}C^{2/(d+\#I^c)}\log n),
\end{aligned}$$

so convergence of the series in (20) for $C$ large enough ($C = \tilde{L}_d = \tilde{L}_d D$) follows from Lemma 13. This completes the proof of Step 5.

To prove sufficiency of (12), by Corollary 1 it is enough to show convergence of the series

$$\sum_n \mathbb{P}\biggl(\biggl|\sum_{|\mathrm{i}|\le 2^n} h(\mathbf{X}_{\mathrm{i}}^{\mathrm{dec}})\biggr| \ge C2^{nd/2}\log^{d/2} n\biggr) \tag{27}$$

for $C = L_d D$. To this end for each $n$ we simply decompose $\Sigma$ into five disjoint sets $A_n^i$, $i = 1,\dots,5$, with $A_n^i$ being a set of the form defined at the $i$th step above (which clearly can be done as the union of the sets from Steps 1–5 is the whole $\Sigma$). For $C = L_d D$, from the triangle inequality and Steps 1–5, we get the convergence of the series

$$\sum_n \mathbb{P}\biggl(\biggl|\sum_{|\mathrm{i}|\le 2^n} \pi_d h(\mathbf{X}_{\mathrm{i}}^{\mathrm{dec}})\biggr| \ge C2^{nd/2}\log^{d/2} n\biggr),$$

which is exactly (27), since by the complete degeneracy $\pi_d h = h$. □

## 7. The undecoupled case.

We are now ready to prove our main result.

PROOF OF THEOREM 1. Sufficiency follows from Corollary 3 and Theorem 5. To prove the necessity assume that (1) holds and observe that from Lemma 7 and Corollary 2, $h$ satisfies the randomized decoupled LIL (8) and thus, by Theorem 5, the growth conditions on functions $\|h\|_{\mathcal{J},u}$ are also satisfied [note that the $\|\cdot\|_{\mathcal{J},u}$ norms of the kernel $h(X_1,\dots,X_d)$ and $\varepsilon_1\cdots\varepsilon_d h(X_1,\dots,X_d)$ are equal]. Thus, the only thing that remains to be proved is the complete degeneracy of $h$. The integrability condition (3) implies that $\mathbb{E}|\pi_d h|^p<\infty$ for all $p<2$ and thus from the Marcinkiewicz type laws of large numbers for completely degenerate $U$-statistics by Giné and Zinn [6] it follows that

$$\frac{1}{n^{d/p}}\sum_{\mathrm{i}\in I_n^d}\pi_d h(\mathbf{X}_{\mathrm{i}})\to 0 \qquad \text{a.s.}$$

as $n\to\infty$. Moreover, from the LIL, we have also

$$\frac{1}{n^{d/p}}\sum_{\mathrm{i}\in I_n^d}h(\mathbf{X}_{\mathrm{i}})\to 0 \qquad \text{a.s.}$$

Let us notice that by Hoeffding's decomposition (Lemma 2),

$$\begin{aligned}&\sum_{\mathrm{i}\in I_n^d}(h(\mathbf{X}_{\mathrm{i}})-\pi_d h(\mathbf{X}_{\mathrm{i}}))\\ &\quad=\sum_{k=0}^{d-1}\binom{d}{k}\cdot\frac{(n-k)!}{n!}\cdot\frac{n!}{(n-d)!}\sum_{\substack{i_1,\dots,i_k\le n\\ i_j\ne i_l\text{ for } j\ne l}}\pi_k h(X_{i_1},\dots,X_{i_k})\\ &\quad=(n-d+1)\sum_{\substack{i_1,\dots,i_{d-1}\le n\\ i_j\ne i_l\text{ for } j\ne l}}g(X_{i_1},\dots,X_{i_{d-1}}),\end{aligned} \tag{28}$$

where

$$g(x_1,\dots,x_{d-1})=\frac{1}{(d-1)!}\sum_{\sigma}\tilde g(x_{\sigma(1)},\dots,x_{\sigma(d-1)}),$$

where the sum is over all permutations of $I_{d-1}$ and

$$\tilde g(x_1,\dots,x_{d-1})=\sum_{k=0}^{d-1}\binom{d}{k}\pi_k h(x_1,\dots,x_k).$$

We thus obtain

$$\frac{n-d+1}{n^{d/p}}\left|\sum_{\substack{i_1,\dots,i_{d-1}\le n\\ i_j\ne i_l\text{ for j}\ne\text{l}}}g(X_{i_1},\dots,X_{i_{d-1}})\right|\to 0 \qquad \text{a.s.}$$

Therefore

$$\frac{1}{n^{d/p-1}}\left|\sum_{\substack{i_1,\dots,i_{d-1}\le n\\ i_j\ne i_l\text{ for j}\ne\text{l}}}g(\mathbf{X}_{\mathrm{i}})\right|$$

is stochastically bounded. Putting $p = 2d/(d+1)$ we obtain the CLT normalization for $U$-statistics of order $d-1$ (see for instance [3], Theorem 4.2.4) and from the results by Giné and Zinn ([7], Theorem 1, or [3], Theorem 4.2.6) we get that $g$ is canonical and $\mathbb{E}g^2 < \infty$. Now the CLT for canonical $U$-statistics yields that

$$g(X_1, \ldots, X_{d-1}) = 0 \qquad \text{a.s.}$$

and (28) for $n = d$ gives $h = \pi_d h$, which proves the complete degeneracy of $h$. $\square$

## 8. Final remarks.

Remark. In Theorem 5 the necessary and sufficient conditions for the decoupled LIL were found, under an additional assumption that the kernel is canonical. We would like to remark that the canonicity actually follows from the decoupled LIL, similarly as in the proof of Theorem 5. The proof would however require developing "a decoupled counterpart" of all the limit theorems for $U$-statistics (like CLT and Marcinkiewicz LLN), which would make it quite lengthy and would not involve genuinely new ideas.

*The cluster set.* When $\mathbb{E}h^2 < \infty$, the limit set in the LIL (1) is almost surely equal to

$$\{\mathbb{E}h(X_1, \ldots, X_d)f(X_1)\cdots f(X_d) : \mathbb{E}f^2(X_1) \le 1\}$$

as is proven in [2]. In general this is not the case. For $d = 2$ it is known that the cluster set is an interval [8], whose end-points are known in some special cases [13]. In these special cases, the lim sup turns out to be a relatively complicated function of the "deterministic" lim sup's appearing in the nasc's conditions. It is natural to conjecture that in general the lim sup is also a function of these quantities.

Now we would like to state the following.

Theorem 6. *The cluster set in the LIL (1) is an interval.*

Proof. It is enough to show that

$$\limsup_{n\to\infty} \left| \frac{\sum_{\mathbf{i}\in I_n^d} h(\mathbf{X}_\mathbf{i})}{n^{d/2}\log\log^{d/2} n} - \frac{\sum_{\mathbf{i}\in I_{n-1}^d} h(\mathbf{X}_\mathbf{i})}{(n-1)^{d/2}\log\log^{d/2}(n-1)} \right| = 0 \qquad \text{a.s.},$$

which will follow if we prove that

$$\limsup_{n\to\infty} \frac{1}{n^{d/2}\log\log^{d/2} n} \left| \sum_{\mathbf{i}\in I_n^d, i_d = n} h(\mathbf{X}_\mathbf{i}) \right| = 0 \qquad \text{a.s.}$$

We can reduce the last statement to

$$\sum_n \sum_{2^{n-1}<k\le 2^n} \mathbb{P}\Big(\Big|\sum_{\mathrm{i}\in I_k^d, i_d=k} h(\mathbf{X}_\mathrm{i})\Big| > \delta 2^{nd/2}\log^{d/2} n\Big) < \infty \tag{29}$$

for all $\delta > 0$. Let $\bar{\pi}_{d-1}$ stand for the Hoeffding projection with respect to the first $d-1$ variables only. Then, the complete degeneracy of $h$, gives $\bar{\pi}_{d-1}h = h$, thus to get (29) it suffices to prove that

$$\sum_n \sum_{2^{n-1}<k\le 2^n} \mathbb{P}\Big(\Big|\sum_{\mathrm{i}\in I_k^d, i_d=k} \bar{\pi}_{d-1}h(\mathbf{X}_\mathrm{i})\Big| > \delta 2^{nd/2}\log^{d/2} n\Big) < \infty.$$

We will now proceed similarly as in the first steps of the proof of Theorem 5, that is we will prove the above convergence with $h$ replaced by $h_n = h\mathbf{1}_{A_n}$ for suitable sets $A_n$.

*Step* 1.

$$A_n \subseteq \{x \in \Sigma^d : h^2(x) \ge 2^{nd}\log^d n\}.$$

Since for $2^{n-1} < k \le 2^n$, $\#\{\mathrm{i} \in I_k^d : i_d = k\} \le 2^{n(d-1)}$ we can use the Chebyshev inequality, exactly as in the first step of the proof of Theorem 5.

*Step* 2.

$$A_n \subseteq \{x : h^2(x) \le 2^{nd}\log^d n, \exists_{I\subseteq I_{d-1}, I\neq\emptyset}\ \mathbb{E}_I(h^2 \wedge 2^{2nd}) \ge 2^{\#I^c n}\log^d n\}.$$

Note that by the decoupling inequalities for the moments of $U$-statistics (see, e.g., [3], Theorem 3.1.1) and Lemma 1 applied conditionally on $X_k$, we have

$$\begin{aligned}\mathbb{E}\Big|\sum_{\mathrm{i}\in I_k^d, i_d=k} \bar{\pi}_{d-1}h(\mathbf{X}_\mathrm{i})\Big| &\le L_d\mathbb{E}\Big|\sum_{\mathrm{i}\in I_k^d, i_d=k} \bar{\pi}_{d-1}h(\mathbf{X}_\mathrm{i}^{\mathrm{dec}})\Big| \\ &\le 2^{d-1}L_d\mathbb{E}\Big|\sum_{\mathrm{i}\in I_k^d, i_d=k} \epsilon_\mathrm{i}^{\mathrm{dec}} h(\mathbf{X}_\mathrm{i}^{\mathrm{dec}})\Big|.\end{aligned}$$

Therefore if we define the sets $A_n(I)$ and $A_{n,l}(I)$ (for $I \subseteq I_{d-1}, I \neq \emptyset$) like in Step 2 of the proof of Theorem 5, it is enough to prove

$$\sum_n \sum_{2^{n-1}<k\le 2^n} \frac{1}{2^{nd/2}\log^{d/2} n}\sum_{l=0}^{\infty}\mathbb{E}\Big|\sum_{\mathrm{i}\in I_k^d, i_d=k} \epsilon_\mathrm{i}^{\mathrm{dec}} h_{n,l}(\mathbf{X}_\mathrm{i})\Big| < \infty,$$

where for fixed $I$ the function $h_{n,l}$ are defined as in the proof of Theorem 5. But for each $2^{n-1} < k \le 2^n, l$ we have by a similar computation as

there

$$\mathbb{E}\left|\sum_{\mathrm{i}\in I_k^d, i_d=k} \epsilon_{\mathrm{i}}^{\mathrm{dec}} h_{n,l}(\mathbf{X}_{\mathrm{i}})\right| \le [2^{(\#I^c+d/2-1)n+l+1}\log^{d/2} n] \\ \times \mathbb{P}_{I^c}(\mathbb{E}_I(h^2\wedge 2^{2nd}) \ge 2^{2l+\#I^c n}\log^d n).$$

Thus

$$\sum_{2^{n-1}<k\le 2^n} \mathbb{E}\left|\sum_{\mathrm{i}\in I_k^d, i_d=k} \epsilon_{\mathrm{i}}^{\mathrm{dec}} h_{n,l}(\mathbf{X}_{\mathrm{i}})\right| \\ \le [2^{(\#I^c+d/2)n+l+1}\log^{d/2} n] \\ \times \mathbb{P}_{I^c}(\mathbb{E}_I(h^2\wedge 2^{2nd}) \ge 2^{2l+\#I^c n}\log^d n)$$

and we can finish this step just as Step 2 in the proof of Theorem 5.

*Step* 3.
$A_n \subseteq \{x : h^2(x) \le 2^{nd}\log^d n, \forall_{I\subseteq I_{d-1}, I\neq\varnothing}\ \mathbb{E}_I(h^2\wedge 2^{2nd}) \le 2^{\#I^c n}\log^d n\}$.

Using the same arguments as above and the Khintchine inequality for Rademacher chaoses we obtain

$$\mathbb{E}\left|\sum_{\mathrm{i}\in I_k^d, i_d=k} \bar{\pi}_{d-1} h(\mathbf{X}_{\mathrm{i}})\right|^4 \le L_d \mathbb{E}\left|\sum_{\mathrm{i}\in I_k^d, i_d=k} \bar{\pi}_{d-1} h(\mathbf{X}_{\mathrm{i}}^{\mathrm{dec}})\right|^4 \\ \le 2^{4(d-1)} L_d \mathbb{E}\left|\sum_{\mathrm{i}\in I_k^d, i_d=k} \epsilon_{\mathrm{i}}^{\mathrm{dec}} h(\mathbf{X}_{\mathrm{i}}^{\mathrm{dec}})\right|^4 \\ \le \tilde{L}_d \mathbb{E}\left(\sum_{|\mathrm{i}|\le k, i_d=k} h^2(\mathbf{X}_{\mathrm{i}}^{\mathrm{dec}})\right)^2,$$

where in the last inequality we have added the diagonal summands just to make the proof more similar to the analogous step (Step 3) in the proof of Theorem 5. Therefore, it suffices to prove

$$\sum_n \frac{2^n \mathbb{E}(\sum_{|\mathrm{i}|\le 2^n, i_d=2^n} h_n^2(\mathbf{X}_{\mathrm{i}}^{\mathrm{dec}}))^2}{2^{2nd}\log^{2d} n} < \infty.$$

But again this can be done just as in Step 3 in the proof of Theorem 5, by considering just the cases (a), (b) there. The case (c) (which made all the consequent work in the proof of Theorem 5 necessary) cannot appear here because the index $i_d$ is fixed. The proof of the theorem may be thus complete just as for Theorem 5, by splitting $\Sigma^d$ into 3 parts (for each $n$), corresponding to Steps 1–3 above. □

Institute of Mathematics
Polish Academy of Sciences
Śniadeckich 8
P.O.Box 21
00-956 Warszawa 10
Poland
E-mail: R.Adamczak@impan.gov.pl

Institute of Mathematics
Warsaw University
Banacha 2
02-097 Warszawa
Poland
E-mail: rlatala@mimuw.edu.pl